\documentclass[12pt]{article}
\usepackage{amsfonts}
\usepackage{amsthm}
\usepackage{mathrsfs}
\usepackage{amsmath}
\usepackage{amssymb}
\usepackage{xcolor}
\usepackage{color}
\usepackage{titlesec}
\usepackage{subeqnarray}
\usepackage[T1]{fontenc}

\newcommand\R{{\mathbb R}}
\newcommand\T{{\mathbb T}}
\newcommand\C{{\mathbb C}}

\def\AA{{\mathcal A}}
\def\BB{{\mathcal B}}

\def\DD{{\mathcal D}}

\def\MM{{\mathcal M}}
\def\NN{{\mathcal N}}
\def\OO{{\mathcal O}}
\def\PP{{\mathcal P}}
\def\QQ{{\mathcal Q}}

\def\UU{{\mathcal U}}
\def\VV{{\mathcal V}}

\def\XX{{\mathcal X}}

\def\ZZ{{\mathcal Z}}

\def\BBB{{\mathscr B}}

\def\LLL{{\mathscr L}}

\def\PPP{{\mathscr P}}

\def\Ppp{{\mathbf P}}

\newtheorem{theorem}{Theorem}[section]

\newtheorem{lemma}[theorem]{Lemma}

\newtheorem{remark}[theorem]{Remark}

\def\Nt{|\hskip-0.04cm|\hskip-0.04cm|}

\def\QEDopen{{\setlength{\fboxsep}{0pt}\setlength{\fboxrule}{0.2pt}\fbox{\rule[0pt]{0pt}{1.3ex}\rule[0pt]{1.3ex}{0pt}}}}
\def\QED{\QEDopen}
\def\endproof{\hspace*{\fill}~\QED\par\endtrivlist\unskip}

\def\eps{{\varepsilon}}

\makeatletter 
\@addtoreset{equation}{section}
\makeatother  

\newcommand{\beqn}{\begin{equation}}
\newcommand{\eeqn}{\end{equation}}
\newcommand{\bear}{\begin{eqnarray}}
\newcommand{\eear}{\end{eqnarray}}
\newcommand{\bean}{\begin{eqnarray*}}
\newcommand{\eean}{\end{eqnarray*}}

\newcommand{\Black}{\color{black}}

\newcommand{\Blue}{\color{black}}

\def\signsm{\medskip\noindent
St\'ephane Mischler.  Universit\'e Paris-Dauphine,   PSL Research University,
CNRS, UMR [7534], CEREMADE, 
Place du Mar\'echal de Lattre de Tassigny
75775 Paris Cedex 16, 
France. Email: mischler@ceremade.dauphine.fr}

\def\signqw{\medskip\noindent
Qilong Weng.  Universit\'e Paris-Dauphine, PSL Research University,
CNRS, UMR [7534], CEREMADE, 
Place du Mar\'echal de Lattre de Tassigny
75775 Paris Cedex 16, 
France. Email: weng@ceremade.dauphine.fr}

\begin{document}

%


\title{Relaxation in time elapsed neuron network
models in the weak connectivity regime}

\author{S. MISCHLER \& Q. WENG}

\date{}

\maketitle

\begin{abstract}
In order to describe the firing activity of a homogenous assembly of
neurons, we consider time elapsed models, which give
mathematical descriptions of the probability density of neurons
structured by the distribution of times elapsed since the last
discharge. Under general assumption on the firing rate and the delay distribution, 
we prove the uniqueness of the steady
state and its nonlinear exponential stability  in the
weak connectivity regime. 
In other words, total asynchronous firing of neurons appears asymptotically in large time. 
The result generalizes
some similar results obtained in \cite{PPS1,PPS2} in the case without delay. Our approach uses the spectral analysis theory for semigroups
in Banach spaces developed recently by the first author and collaborators. 
\end{abstract}

\begin{center} {\bf Version of \today}
\end{center}

\smallskip
\noindent \textbf{Keywords.} Neuron networks, time elapsed dynamics,
semigroup spectral analysis, weak connectivity, long time
asymptotic.

 \smallskip 
\noindent \textbf{AMS Subject Classification.} 35B40, 35F15, 35F20, 92B20 

\tableofcontents

\section{Introduction}
\label{sec:Intro}

In nervous systems, neuronal circuits carry out tasks of information
transmission and processing. Many neurons generate trains of
stereotyped electrical pulses in response to incoming stimulations.
Following each discharge, the neuron undergoes a period of
refractoriness during which it is less responsive to inputs,
before recovering its excitability \cite{PPS2}. The main carrier
of information is the discharge times or some statistics of the
discharge times. 
In this work, we consider  a simple neuronal model which neglects the mechanisms underlying spike generation and 
focusses on describing the neuronal dynamics in terms of discharge times. More precisely, we consider a model
which has been introduced and studied in \cite{GKbook,PPS1,PPS2} and which 
describes the post-discharge recovery of neuronal
membranes through an instantaneous firing rate that depends on the
time elapsed since the last discharge and the inputs by neurons.
We refer to these papers for biologic motivation
and discussions. We also refer to \cite{MR3311484,FL*,RT*,Q*} where
these models (or similar ones) are obtained as a mean field limit
of finite number of neuron network models.  


The neuronal network is described here by the density number of neurons
$f=f(t,x)\geq0$ which at time $t \ge 0$ are in the state $x \ge 0$. The state of a neuron is a local time (or internal clock)  
which corresponds to the elapsed time since the last
discharge. The dynamic of the neuron network 
is given by the following nonlinear time elapsed (or of age structured type) evolution  equation
\begin{subequations}\label{ASM}
 \begin{align}
   &\partial_t f=-\partial_x f-a(x,\eps \,  m(t))f=:\mathcal{L}_{\eps m(t)}f,\\
   &f(t,0)=p(t), \ \ f(0,x)=f_0(x).
 \end{align}
\end{subequations}
Here $a(x,\eps \, \mu)\geq0$ represents the firing rate of a neuron in
the state $x$  for a network activity $ \mu\geq0$  and a network connectivity parameter  $\eps \ge 0$. 
The function $p(t)$ represents the total density of neurons
which undergo a discharge at time $t$ and is defined through
\beqn\label{eq:ASM2}
p(t):=\mathcal{P}[f(t); m(t)], 
\eeqn
 where 
\beqn\label{eq:ASM3}
 \PP [g,\mu] = \PP_\eps[g,\mu] :=  \int_0^\infty
    a(x,\eps \mu) g(x)\mathrm{d}x.
\eeqn
The function $ m(t)$ represents the 
network
activity at time $t\geq0$ resulting from earlier discharges and
is defined by
\beqn\label{eq:delaym=pb}
    m(t):=\int_0^\infty p(t-y)b(\mathrm{d}y),
\eeqn
where the delay distribution $b$ is a probability measure which takes into account the persistence of the electric activity in the network resulting from discharges (synaptic integration). In the sequel, we will consider the two following situations :

$\bullet$ The {\it case without delay},  when $b = \delta_0$ and then $m(t) = p(t)$.

$\bullet$ The {\it case with delay},  when $b$ is a smooth function.

\smallskip
We observe that in both cases, the solution $f$ of the time elapsed equation
\eqref{ASM}--\eqref{eq:ASM3} satisfies
$$
\frac{\mathrm{d}}{\mathrm{d}t}\int_0^{\infty}f(t,x)\mathrm{d}x=f(t,0)-
\int_0^{\infty}a(x,\eps m(t))f(t,x)\mathrm{d}x=0. 
$$
As a consequence, the total density number of neurons (also called {\it mass} in the sequel) is
conserved and we can normalize that mass to be $1$. In other words, we may always assume  
$$
\langle f(t,.) \rangle = \langle f_0 \rangle = 1,
\quad \forall t\geq0, \quad \langle g \rangle := \int_0^{\infty} g(x)\mathrm{d}x.
$$

\smallskip
A (normalized) steady state for the  time elapsed evolution system of equations \eqref{ASM}--\eqref{eq:ASM3} is a couple $(F_\eps,M_\eps)$ of a density number of neurons
$F_\eps = F_\eps(x) \geq0$ and a network activity $M_\eps \ge 0$ such that 
\begin{subequations}\label{eq:StSt}
 \begin{align} \label{eq:StSt1}
   &0=-\partial_x F_\eps - a(x,\eps \,  M_\eps) F_\eps= \mathcal{L}_{\eps M_\eps} F_\eps,\\
   &F_\eps(0)= M_\eps, \quad \langle F_\eps \rangle = 1. 
 \end{align}
\end{subequations}
It is worth emphasizing that for a steady state the associated network
activity and discharge activity are two equal constants because of the
normalization of the delay distribution, i.e. $\langle b\rangle=1$. 

\smallskip
In equations \eqref{ASM}--\eqref{eq:ASM3} and \eqref{eq:StSt}, the connectivity parameter $\eps \ge 0$ corresponds to the strength of the influence 
of   the neuronal network activity on  each neuron through the functions $m(t)$ and $f(t,x)$ respectively. In the limit case $\eps = 0$, equation 
 \eqref{ASM}--\eqref{eq:ASM3} is linear which means that  each neuron evolves accordingly to its own dynamic. In the other hand, when $\eps > 0$, equation 
  \eqref{ASM}--\eqref{eq:ASM3} is nonlinear and the dynamic of any given neuron is affected by the state (or the past states in the case of the model with delay) of all the other neurons through the global  activity 
  of the neuronal network. Finally, the weak connectivity regime, about which we are mainly concerned in the present paper, 
  corresponds to a range of connectivity parameter $\eps \in (0,\eps_0]$, $\eps_0 > 0$ small enough, such that the nonlinearity of equations \eqref{ASM}--\eqref{eq:ASM3} and \eqref{eq:StSt} is 
  not too strong. 

\smallskip
Our main purpose in this paper is to prove that solutions to the  time elapsed evolution  equation \eqref{ASM}--\eqref{eq:ASM3} converge to a stationary state under a
weak connectivity assumption. Before stating that result, let us present the precise mathematical assumptions we will need on the firing rate $a$ and 
on the delay distribution $b$. 

\smallskip
We make the physically reasonable assumption
\beqn\label{hyp:a1}
\partial_x a\geq0, \ \ a'=\partial_\mu a\geq0,
\eeqn
\beqn\label{hyp:a2}
0<a_0:=\lim_{x\rightarrow\infty}a(x,0)\leq\lim_{x,\mu\rightarrow\infty}a(x,\mu)=:a_1<\infty,
\eeqn
as well as the smoothness assumption
\begin{equation}\label{hyp:a3}
a \in W^{2,\infty}(\R^2_+). 
\end{equation}
In the delay case, we assume that $b(dy) = b(y) \, dy$ satisfies the exponential bound and smoothness condition
\begin{equation}\label{hyp:del}
\exists\delta>0, \quad\int_0^\infty e^{\delta
y}\, (b(y) +  |b'(y)|) \, \mathrm{d}y <\infty.
\end{equation}

 We begin by stating our main result about the stationary problem~\eqref{eq:StSt}.

\begin{theorem}\label{th:SS}
Assume \eqref{hyp:a1}-\eqref{hyp:a2}-\eqref{hyp:a3}. \Blue For any $\eps \ge 0$, there exists at least one solution 
$(F_\eps,M_\eps)\in W^{1,\infty}(\R_+) \times\mathbb{R}_+$
to the stationary problem \eqref{eq:StSt} such that 
\beqn\label{eq:F&FprimBd}
0 \le F_\eps (x)  \le C \, e^{-{a_0 \over 2} x}, \quad  |F_\eps'(x)| \le C \, e^{-{a_0 \over 2} x}, \quad \forall \, x \ge 0,
\eeqn
for a constant $C \in (0,\infty)$. \Black
Moreover, there exists $\eps_0 > 0$, small enough, such that 
the above solution is unique for any
$\eps\in[0,\eps_0)$.
\end{theorem}

\smallskip 
For a given initial datum $0 \le f_0 \in L^1(\R_+)$, we say that a function $f$ is a weak (positive and mass conserving)  solution to \eqref{ASM}--\eqref{eq:ASM3} if 
$$
0 \le f \in C([0,\infty); L^1(\R)), \quad \langle f(t) \rangle =  \langle f_0 \rangle, \,\, \forall \, t \ge 0, 
$$
and $f$ satisfies  \eqref{ASM} in the distributional sense $\DD'([0,\infty) \times [0,\infty))$  for some functions $m,p \in C([0,\infty))$ which fulfilled the constrains \eqref{eq:ASM2} and
\eqref{eq:ASM3}. 

\smallskip
Under the above assumptions, existence and uniqueness of  weak solutions  have been established in  \cite[Theorem 1.1]{W*}. 
The main concern of the present work is the following long-time asymptotic result on the solutions. 

\begin{theorem}\label{th:MR} We assume that the firing rate $a$ satisfies \eqref{hyp:a1},  \eqref{hyp:a2} and  \eqref{hyp:a3}. We also assume
that the  delay distribution  $b$  satisfies $b = \delta_0$ or  \eqref{hyp:del}. 
There exists $\eps_0 > 0$, small enough, and there exist some constants $\alpha < 0$, $C \ge 1$  and  $\eta > 0$  such that for any connectivity parameter $\eps \in (0,\eps_0)$  and any initial datum $0 \le f_0\in L^1$ with mass $1$  and such that $\| f_0 - F_\eps \|_{L^1} \le \eta/\eps$, the unique 
solution $f$ to the evolution equation \eqref{ASM}--\eqref{eq:ASM3} 
satisfies 
$$
\|f(t,.)-F_\eps \|_{L^1}\leq C  e^{\alpha t}, \qquad \forall \, t \ge 0 .
$$
\end{theorem}

In other words, in that weak connectivity regime,
 we prove that the total asynchronous firing of neurons appears exponentially fast in the large time asymptotic. 
Theorem~\ref{th:MR} extends to firing rates $a$ satisfying  \eqref{hyp:a1}--\eqref{hyp:a3} some similar exponential stability results obtained in \cite{PPS1,PPS2}
in the case without delay and for a  {\it step function} firing rate $a$ given by
\beqn\label{eq:StepFctStructure}
a(x,\mu) = {\bf 1}_{x > \sigma(\mu)}, \quad \sigma,\sigma^{-1}\in W^{1,\infty}(\R_+), \quad \sigma' \le 0.
\eeqn
It is worth mentioning that the above  firing rate does not fall in the class of rates considered in the present paper because condition \eqref{hyp:a3} is not met. 
On the other hand, we are able to tackle the case without and with delay in the same time, what  it was not the case in \cite{PPS1,PPS2}. In the delay case,  stability results were established in \cite{PPS1}, but not exponential stability. 

\smallskip
{Our proof follows a strategy of ``perturbation of semigroup'' initiated in \cite{MMcmp} for studying long time convergence to the equilibrium for the homogeneous inelastic Boltzmann equation and used recently in \cite{MQT} for a neuron network equation based on a brownian (hypoelliptic) perturbation of the well-known FitzHugh-Nagumo dynamic.  More precisely, we introduce the linearized equation for the variation functions $(g,n,q) = (f,m,p) - (F_\eps,M_\eps,M_\eps)$ around a stationary state $(F_\eps,M_\eps,M_\eps)$, which writes
\begin{subequations}\label{eq:ASMlin1}
 \begin{align}
   &\partial_t g=-\partial_x g- a(x,\eps \, M_\eps) g - n(t) \, \eps (\partial_\mu a)(x,\eps M_\eps) F_\eps,  \\
   &g(t,0)=q(t), \ \ g(0,x)= g_0(x),
    \end{align}
\end{subequations}
with  
\beqn\label{eq;ASMlin2}
q(t) =   \int_0^\infty a(x,\eps \, M_\eps) g \,\mathrm{d}x +  n(t) \, \eps \int_0^\infty  (\partial_\mu a)(x,\eps M_\eps) F_\eps \,  \mathrm{d}x
\eeqn
and
\beqn\label{eq;ASMlin3}
  n(t):=\int_0^\infty q(t-y)b(\mathrm{d}y). 
\eeqn
We associate to that linear evolution equation a generator $\Lambda_\eps$ (which acts on an appropriate space to be specified in the
two cases without and with delay) and its semigroup
$S_{\Lambda_\eps}$. It turns out that we may split the operator $\Lambda_\eps$ as
$$
\Lambda_\eps=\AA_\eps+\BB_\eps, 
$$
for some $\alpha$-hypodissipative operator $\BB_\eps$, $\alpha<0$, and some bounded and  $\BB_\eps$-power regular operator
$\AA_\eps$ as defined in \cite{Voigt80,GMM*,MS,MS-erratum,Mbook*}. In particular, adapted versions of the Spectral Mapping Theorem in \cite{MS,MS-erratum,Mbook*} and the   Weyl's Theorem in
 \cite{Voigt80,GMM*,MS,MS-erratum,Mbook*} imply that the semigroup
 $S_{\Lambda_\eps}$ has a finite dimensional dominant part. Moreover, 
 in the limit case when $\eps=0$, the term $n(t)$ disappears
 from equation \eqref{eq:ASMlin1} and the resulting semigroup $S_{\Lambda_0}$ becomes positive. That allows us to use the Krein-Rutman Theorem established in \cite{MS,MS-erratum,Mbook*} in order to get that the stationary state $(F_0,M_0,M_0)$ is unique and exponentially stable. Using next a perturbative argument developed in \cite{MMcmp,Tristani,Mbook*,MT*}, we get that the unique stationary state $(F_\eps,M_\eps,M_\eps)$ is also exponentially stable in the weak connectivity regime. 
We conclude the proof of Theorem~\ref{th:MR} by a somewhat classical nonlinear exponential stability argument.  

\smallskip

The same strategy applies to the case without delay and with delay. In both cases, the boundary condition in the age structure equation is treated as a source term and, in the delay case, the delay equation~\eqref{eq:delaym=pb} is replaced  by a simple age equation on an auxiliary function, so that the resulting linearized equation writes as an autonomous system of two PDEs and falls in the classical framework of linear evolution equations generating a semigroup. 

\smallskip
Our approach is thus quite different from the usual way to deal with delay equations, as introduced by I.~Fredholm~\cite{MR1554993} and V.~Voltera~\cite{MR0100765}, which consists in using the specific framework of   ``fading memory space", which goes back at least to Coleman \& Mizel   \cite{MR0210343},  or the theory of ``abstract algebraic-delay differential systems" developed by O. Diekmann and co-authors \cite{MR1345150}.
 
\smallskip
Our approach is also different from the previous works \cite{PPS1,PPS2,PPS3} where the asymptotic stability analysis were performed by taking advantage of the step function structure \eqref{eq:StepFctStructure} of the firing rate. That one makes possible to explicitly exhibit a suitable norm (related to the $W_1$  Monge-Kantorovich-Wasserstein optimal transport distance) such that some related linear age structure operator is dissipative. The present method is based on a more abstract approach but in the other hand it is somewhat more flexible because if does note require to explicitly exhibit a norm for which the underlying linear(ized) operators are dissipative. 
In particular, we {\Blue hope} that our strategy can be adapted to the large connectivity regime, to the step function firing rate  \eqref{eq:StepFctStructure} as well as to 
models including fragmentation term to describe neuronal networks with adaptation  and fatigue, and thus generalize to the case  with a delay term all  the stability results established in \cite{PPS1,PPS2,PPS3} in the case without delay. 

\medskip
Let us end the introduction by describing the plan of the paper. 
In Section~\ref{sec:WithoutDelay}, we introduce the strategy, we prove the stationary state result and we establish Theorem~\ref{th:MR} in the case without delay.
In Section~\ref{sec:WithDelay}, we  establish Theorem~\ref{th:MR} in the case with  delay. 
As we mentioned above, the  strategy of proof for the case without and with delay is rather the same. For pedagogical reason, we start presenting the method on the simplest {\it ``without delay case"}  in 
Section~\ref{sec:WithoutDelay}, where we prove the  stationary problem result Theorem~\ref{th:SS} as well as Theorem~\ref{th:MR} in that case.
Next, in Section~\ref{sec:WithDelay}, we only explain how the proof must be modified in order to treat the more complicated {\it``with delay case"} and thus 
 establish Theorem~\ref{th:MR} in all generality. 

\medskip
\paragraph{\bf Acknowledgments.}  The research leading to this paper was (partially) funded by the French "ANR blanche" project Kibord: ANR-13-BS01-0004.  

\section{Case without delay}
\label{sec:WithoutDelay}

The present section is devoted to the proof of our main result Theorem~\ref{th:MR} in the  without delay case.

\subsection{The stationary problem}

We first deal with the stationary problem and we prove the existence of steady state as well as its uniqueness in the small connectivity regime.

\proof[Proof of Theorem \ref{th:SS}]  {\sl Step 1. We prove the existence of a solution.}
We set
$$
A(x,m):=\int_0^x a(y,m)\mathrm{d}y, \quad \forall \ x,m \ge 0. 
$$
 For any
$m\geq0$, we can solve the equation \eqref{eq:StSt1}, by writing 
\beqn\label{eq:Fepsm=}
F_{\eps,m}(x):=T_m e^{-A(x,\eps m)},
\eeqn
where $T_m \ge 0$ is chosen in order that  $F_{\eps,m}$ satisfies the mass normalized condition, 
namely 
$$
T_m^{-1}=\int_0^{\infty}e^{-A(x,\eps m)}\mathrm{d}x.
$$
In order to conclude the existence of a solution, we just have to find a real number $m = M_\eps$ such that 
$m = F_{\eps,m}(0) = T_m$. Equivalently, we need to find $M_\eps \ge 0$ such that
\begin{equation}\label{mnc}
\Phi(\eps,M_\eps)=1,
\end{equation}
where
$$
 \Phi(\eps,m)=m
T_{m}^{-1} := m\int_0^{\infty}e^{-A(x,\eps m)}\mathrm{d}x.
$$

From the assumption \eqref{hyp:a2} of $a $, there
exists $x_0\in[0,\infty)$ such that $a(x,\mu)\geq\frac{a_0}{2}$, for
any $x\geq x_0$, $\mu \ge 0$, and therefore
\begin{equation}\label{EoA}
    \frac{a_0}{2}(x-x_0)_+\leq A(x,\mu)\leq a_1 x, \ \ \forall \, x\geq0, \,\, \forall \, \mu \ge 0.
\end{equation}

\Blue 
Thanks to the Lebesgue dominated convergence theorem, we deduce that $\Phi(\eps,\cdot)$ is a continuous function. 
Because $\Phi(\eps,0)=0$, $\Phi(\eps,\infty)=\infty$ and thanks to the intermediate value theorem, we conclude to the existence of at least one real number $M_\eps \in (0,\infty)$
 such that \eqref{mnc} holds. Estimates \eqref{eq:F&FprimBd} immediately follow from the identity \eqref{eq:Fepsm=} and the estimate \eqref{EoA}. \Black
%

\smallskip\noindent   {\sl Step 2. We prove the uniqueness of the solution in the weak connectivity regime.}
 Obviously, there exists a unique $M_0 := (\int_0^{\infty}e^{-A(x,0)}\mathrm{d}x)^{-1}\in(0,\infty)$
such that $\Phi(0,M_0)=1$. Moreover, we compute
$$
\frac{\partial}{\partial
m}\Phi(\eps,m)=\int_0^{\infty}e^{-A(x,\eps
m)} (1 -m\eps {\partial A \over \partial m} (x,\eps m))  \mathrm{d}x,
$$
which is continuous as a function of the two variables because of  \eqref{hyp:a3}. We then easily obtain that $\Phi \in C^1$. Since moreover
$$
\frac{\partial}{\partial
m}\Phi(\eps,m)|_{\eps=0}=\int_0^{\infty}e^{-A(x,0)}\mathrm{d}x>0,
$$
the implicit function theorem implies that there exists
$\eps_0>0$, small enough, such that  the equation \eqref{mnc}  has a unique solution for any $\eps\in[0,\eps_0)$. 
\endproof

\begin{remark}
In the above proof, we do not need \eqref{hyp:a3} but only the  
weaker smoothness assumption that $A$ and $\partial_m A$ are continuous.
\end{remark}

\subsection{Linearized equation and structure of the spectrum}  \label{subsec:WithoutD:LinearizedEq&structure}
To go one step further, we introduce the linearized equation around the
stationary solution $(F_\eps,M_\eps)$. On the variation $(g,n)$, the linearized equation writes
\bean
   &&\partial_t g+\partial_x g+a_\eps g+a'_\eps F_\eps n(t)=0,\\
   &&g(t,0)=n(t)=\int_0^{\infty}(a_\eps g+a'_\eps F_\eps n(t))\, \mathrm{d}x, 
   \quad g(0,x)=g_0(x),
\eean
with $a_\eps :=a(x,\eps M_\eps)$, $a'_\eps := \eps \, (\partial_\mu a)(x,\eps M_\eps)$. Since there exists $\eps_0 > 0$, small enough, such that 
$$
\forall \, \eps \in (0,\eps_0) \qquad \kappa := \int_0^{\infty} a'_\eps F_\eps \mathrm{d}x < 1,
$$
we may define
\beqn\label{def:Nepsg}
\MM_\eps[g] := (1-\kappa)^{-1} \int_0^\infty a_\eps g\, \mathrm{d}x, 
\eeqn
and the linearized equation is then equivalent to
\bear \label{ASMDL1}
   &&\partial_t g+\partial_x g+a_\eps g+a'_\eps F_\eps \, \MM_\eps[g(t,.)] =0,
   \\
   &&g(t,0)=  \MM_\eps[g(t,.)],\label{ASMDL2}
   \quad g(0,x)=g_0(x). 
\eear
 
\Blue
To the above linear evolution equation, one can classically associate a semigroup $S_{L_\eps}(t)$ acting on $ X := L^1(\R_+)$, with generator
$$
L_\eps g := - \partial_x g - a_\eps g- a'_\eps F_\eps \, \MM_\eps[g] 
$$
and domain
$$
D(L_\eps) := \{ g \in W^{1,1}(\R_+); \,\, g(0) = \MM_\eps[g] \}.
$$

Here we also use another approach by  considering the boundary term as a source term, and then  rewriting the equation as
\beqn\label{eq:ASMDL3}
\partial_t g=\Lambda_\eps g :=-\partial_x g-a_\eps g-a'_\eps F_\eps\MM_\eps[g]+\delta_{x=0}\MM_\eps[g],
\eeqn
acting on the  space  of  bounded Radon measures 
$$
\XX :=M^1(\mathbb{R}_+) = \{ g \in (C_0(\R))'; \,\, \hbox{supp} \, g \subset \R_+ \},
$$
endowed with the weak $*$ topology $\sigma(M^1,C_0)$. Here and below, $C_0(I)$ denotes the space of continuous functions on the closed
interval $I = \R$ or $I = \R_+$ which goes to $0$ at infinity. We then denote by $S_{\Lambda_\eps}(t)$ the semigroup on $\XX$ generated by $\Lambda_\eps$. 

 \smallskip
It is worth emphasizing that for any $g_0 \in X$, the function $g (t) = S_{L_\eps}(t) g_0 \in C([0,\infty); X)$ is a weak solution to equations \eqref{ASMDL1}-\eqref{ASMDL2}, and more precisely clearly satisfies
$$
- \int_0^\infty \varphi(0) \, g_0 \, dx + \int_0^\infty\!\!\! \int_0^\infty g \, \{ - \partial_t \varphi + \Lambda_\eps^* \varphi \} \, dx dt = 0,
$$
for any $\varphi \in C^1_c([0,\infty); C_0(\R_+)) \cap C_c([0,\infty); C^1_0(\R_+))$. Here, we have defined 
$$
\Lambda_\eps^* \psi  :=  \partial_x \psi - a_\eps \, \psi  - a_\eps \, (1-\kappa)^{-1} \, \Bigl[ \int_0^\infty a_\eps'(y) F_\eps(y) \psi(y) \, dy - \psi(0) \Bigl],
$$
for any $\psi \in C^1_0(\R_+)$, the space of $C^1$ functions which goes to $0$ at infinity as well as their first derivative. 
As a consequence, the semigroup $S_{\Lambda_\eps}$ being defined by duality from the semigroup $S_{\Lambda^*_\eps}$, we have $S_{\Lambda_\eps}|_{X} = S_{L_\eps}$. 

\Black
\smallskip
For a generator $L$, we denote by $\Sigma(L)$ its spectrum and by $S_L$ the associated semigroup. 
We refer to the classical textbooks \cite{Kato,Pazy,EngelNagel} for an introduction to the spectral analysis of operators and the semigroup theory.  
 Our next result deals with the structure of the spectral set $\Sigma(\Lambda_\eps)$ of $\Lambda_\eps$ and the splitting structure of the associated semigroup $S_{\Lambda_\eps}$.

\Blue
\Black

\begin{theorem}\label{th:MRe}
Assume \eqref{hyp:a1}-\eqref{hyp:a2}-\eqref{hyp:a3} and define  $a^* := - a_0/2 < 0$.
\Blue The operator $\Lambda_\eps$ is the generator of a   weakly $*$ continuous semigroup $S_{\Lambda_\eps}$ acting on  $\XX$ endowed with the weak $*$ topology $\sigma(M^1,C_0)$. 
Moreover,  there exists a finite rank projector $\Pi_{\Lambda_\eps,a^*}$ which commutes with $S_{\Lambda_\eps}$, an integer $j \ge 0$ and some complex numbers 
$$
\xi_1, ..., \xi_j \in \Delta_{a^*} := \{ z \in \C, \,\, \Re e \, z > a^* \},
$$
such that on $E_1 := \Pi_{\Lambda_\eps,{a^*}} \XX $ the restricted operator satisfies 
\beqn\label{eq:SigmaLambdaE1}
\Sigma(\Lambda_{\eps|E_1} )\cap\Delta_{a^*}=\{\xi_1,...,\xi_j\} 
\eeqn
(with the convention $\Sigma(\Lambda_{\eps|E_1} )\cap\Delta_{a^*} = \emptyset$ when $j=0$) and for any $a > {a^*}$ there exists a constant $C_a$ such that the remainder semigroup satisfies  
\beqn\label{eq:SLambdaepsPiPerp}
\|S_{\Lambda_\eps}(I-\Pi_{\Lambda_\eps,{a^*}})\|_{\mathscr{B}( \XX)}\leq
C_a e^{a t}, \quad \forall \, t \ge 0. 
\eeqn
\end{theorem}

The proof of the result  is a direct consequence of the fact that the  operator $\Lambda_\eps$ splits as $\Lambda_\eps = \AA_\eps + \BB_\eps$ where  $\AA_\eps$ and $\BB_\eps$ are defined on $\XX$ by
\bear\label{def:Aeps}
   && \AA_\eps g:= \gamma_\eps \MM_\eps[g], \quad \gamma_\eps := \delta_0 -a'_\eps F_\eps, \\ \label{def:Beps}
   && \BB_\eps g:=-\partial_x g-a_\eps g,
\eear
for which we can adapt Weyl's Theorem of \cite{Voigt80,GMM*,MS,MS-erratum,Mbook*} and the  Spectral Mapping Theorem of \cite{MS,MS-erratum,Mbook*}.
 The picture may seem not that simple,
because $\Lambda_\eps$ does not generate a strongly continuous
semigroup on $\XX$ and then apparently does not fit the framework developed in
\cite{MS,MS-erratum}. However, we may probably circumvent that issue in the following ways:

\smallskip \Blue
- We may observe that the strong continuity is little used  in  \cite{MS,MS-erratum}  and then the results stated therein extend 
to weakly $*$ continuous semigroups. In other words, on the finite dimensional eigenspace associated to the principal part of the spectrum \eqref{eq:SigmaLambdaE1} continuity and weak $*$ continuity are equivalent, while in the remainder part \eqref{eq:SLambdaepsPiPerp}
we just use a decay bound which does not require the strong continuity property, see \cite[Chapter 1]{Pazy} as well as \cite{W*,Mbook*}
where such a weak $*$ continuous framework is also discussed.

\smallskip
- We may apply the theory developed  in  \cite{MS,MS-erratum} to the adjoint operator $\Lambda^*_\eps$ and the associated semigroup acting on $C_0(\R)$ and then deduce the result by duality.

\smallskip
- We may probably reduce the space $X$ by ``sun duality", see \cite[Chapter II.2.6]{EngelNagel}, apply the theory developed in \cite{MS,MS-erratum} to the resulting  strongly continuous semigroup and  then conclude by a density argument.  

\medskip
We rather follow another strategy, which has the advantage to be more self-contained and pedagogical, by adapting to our context the proofs from \cite{MS,MS-erratum}. 
It is worth emphasizing that we introduce the two spaces $X \subset \XX$ because the splitting $\Lambda_\eps = \AA_\eps + \BB_\eps$ has some clearer meaning in the space $\XX$ while the Banach structure of $X$ is used when we establish the adequate version of Weyl's theorem. 

\Black
\medskip
We recall the definition of hypodissipativity introduced in \cite{GMM*}. 
We say that the closed operator $L$ on a Banach space $E$ with dense domain $D(L)$ is $\alpha$-hypodissipative 
if there exists an equivalent norm $\Nt\cdot\Nt$ on $E$ such that 
$$
   \forall \, f \in D(L), \,\, \exists \, \varphi \in F_{\Nt \cdot \Nt}  (f) \qquad \Re e
  \, \langle \varphi, (L-\alpha) \, f \rangle \le 0,
$$
where, for any $f \in E$, the associated dual set $F_{\Nt \cdot \Nt}  (f)  \subset E'$ is defined by 
$$
F_{\Nt \cdot \Nt} (f):= \{ \varphi \in E'; \,\, \langle  \varphi , f \rangle = \Nt f \Nt_X^2 =\Nt \varphi \Nt_{E'}^2 \}.
$$
We also recall that the generator $L$ of a semigroup  of bounded operators is $\alpha$-hypodissipative if, and only if, there exists a constant $M \ge 1$ such that  the associated semigroup $S_L$  on $E$ satisfies the growth estimate
$$
\| S_L(t) \|_{\BBB(E)} \le M \, e^{\alpha \, t}, \quad \forall \, t \ge 0,
$$
where $\BBB(E)$ denotes the space of linear and bounded operators on $E$. We will sometime abuse by saying that 
$S_L$ is  $\alpha$-hypodissipative  when it satisfies the above growth estimate. We refer to \cite{Pazy,GMM*,MS,Mbook*} for
details. 

\smallskip
We start with the properties of the two auxiliary operators.

\begin{lemma}\label{lem:Lm1}
Assume that $a$ satisfies conditions \eqref{hyp:a1}--\eqref{hyp:a3}. 
For any $\eps \ge 0$, the operators $\AA_\eps$ and
$\BB_\eps$ satisfy : \Blue 
\begin{enumerate}
  \item[(i)]  $\AA_\eps\in\mathscr{B}( W^{-1,1}(\R_+),\XX)$. 
  \item[(ii)] $S_{\BB_\eps}$ is ${a^*}$-hypodissipative in both $X$ and $\XX$.
  \item[(iii)]   The family of operators
  $S_{\BB_\eps}\ast\AA_\eps S_{\BB_\eps}$ satisfies
  $$
  \|(S_{\BB_\eps}\ast\AA_\eps S_{\BB_\eps})(t) \|_{\BBB(\XX,Y)} \le  C_a \, e^{a \, t} , \quad \forall \, a > {a^*}, 
 $$
 for some  some constant $C_a \in (0,\infty)$ and with $Y := BV(\R_+)  \cap L^1_1(\R_+)$.  
\end{enumerate}
Here $BV(\R_+)$ stands for the space of functions with bounded variation and $L^1_1(\R_+)$ stands for the weighted Lebesgue space associated to the
weight function $x \mapsto \langle x \rangle$.  \Black
\end{lemma}

\proof In order to shorten notation we skip the $\eps$ dependency, we write $a (x) = a(x,\eps M_\eps)$, $A(x) = A(x,\eps M_\eps)$, $S_{\Lambda} = S_{\Lambda_\eps}$, 
$S_{\AA} = S_{\AA_\eps}$, $S_{\BB} = S_{\BB_\eps}$ and so on. 

\Blue
\smallskip\noindent{\sl Step 1. Proof of (i). } We have $\AA \in\mathscr{B}(W^{-1,1}(\R_+),\XX)$ from the fact that
$\MM_\eps[\cdot] \in\mathscr{B}(W^{-1,1}(\R_+),\R)$ because $\|a\|_{W^{1,\infty}} < \infty$ as assumed in \eqref{hyp:a3}.
  
\smallskip\noindent{\sl Step 2. Proof of (ii). } 
We use the same notation $\BB$ and $S_\BB$ for the generator and its associated semigroup defined in the spaces $X$ and $\XX$. 
We write $S_{\BB}$ with the explicit formula
\beqn\label{eq:defSBeps1}
  S_{\BB}(t)g(x)=e^{\Blue A(x-t)-A(x)}g(x-t)\mathbf{1}_{x-t\geq0}=:S(t,x).
\eeqn\Blue
 From the  inequality $a (z,0) \ge (3a_0/4) ( 1 - {\bf 1}_{0 \le z \le x_1})$ for some $x_1 \in (0,\infty)$ coming from \eqref{hyp:a2}, we deduce $A(x-t) - A(x) \le 3a_0x_1/4 - (3a_0/4)t $ for any $x \ge t \ge 0$, next
\beqn\label{eq:eA-A}
e^{A(x-t) - A(x)} \le C \, e^{3 \beta t}, \quad \forall \, x \ge t \ge 0, 
\eeqn
with $C :=e^{\frac{3a_0x_1}{4}}>0$, $\beta:=- {a_0}/{4}<0$, and finally
   \begin{eqnarray*}
  \| S_{\BB}(t)g\|_X  
\le Ce^{3\beta t}\|g\|_X, \quad \forall \, t \ge 0, \,\, \forall \, g \in X.
  \end{eqnarray*}
We conclude by observing that $3\beta < {a^*}$ and that the same estimate holds with $X$ replaced by $\XX$ thanks to a (weakly $*$) density argument. 
  
\medskip\noindent{\sl Step 3. Proof of (iii). } For $g \in \XX$ and with the notation of Step (ii), we have 
  $$
 \AA S_{\BB} (t)  g = \gamma_\eps \, N(t), 
 $$
{ where $\gamma_\eps$ is defined in \eqref{def:Aeps}} and
 $$
 N(t) := \MM_\eps[{ S(t,\cdot)}] =(1-\kappa)^{-1} \int_0^\infty a(x) e^{A(x-t)-A(x)}g(x-t)
 \mathbf{1}_{x-t\geq0} \, \mathrm{d}x.
 $$
 It is worth noticing that $N \in C_b(\R_+)$, because $(x,t) \mapsto a(x) e^{A(x-t)-A(x)}$ is a bounded and continuous function
 on the set $\T := \{ (x,t); \, x \ge t \ge 0 \}$. Moreover, as in Step (ii), we have 
\beqn\label{eq:N<ea}
 |N(t)| \le C \, a_1  \int_0^\infty   e^{3\beta t} \,
 |g(x-t)|\mathbf{1}_{x-t\geq0} \, \mathrm{d}x \le C a_1 \, e^{3\beta t} \|g\|_\XX, 
\eeqn
  for any $t \ge 0$. 
 We deduce 
\bear\nonumber
 (S_{\BB}\ast\AA S_{\BB})(t)g(x) 
   &=& \int_0^t (S_{\BB}(s) \gamma_\eps)(x) N(t-s) \,  \mathrm{d}s \\ \label{eq:identitySBASB1} 
   &=& \int_0^t e^{A(x-s)-A(x)} \, \gamma_\eps(x-s)  N(t-s) \mathbf{1}_{x-s\geq0}\,  \mathrm{d}s \qquad\\ \label{eq:identitySBASB2} 
   &=& e^{-A(x)} \, (\nu_\eps * \check N_t) (x), 
\eear
 with $\nu_\eps := \gamma_\eps \, e^A$ and the classical notation $\check N_t(s) = N(t-s)$. 
 Starting from identity \eqref{eq:identitySBASB1} and  denoting $e_{-\beta}(x) := e^{-\beta x}$, we have 
 \begin{eqnarray*}
&&
| (S_{\BB}\ast\AA S_{\BB})(t)g(x) | \, e_{-\beta}(x)
\le e^{-\beta x} \int_0^t C \, e^{3 \beta t} \, |\gamma_\eps(x-s)|   |N(t-s) |  \,  \mathrm{d}s
\\
&&\quad= C \, e^{{a^*} t}   \int_0^t  e^{2\beta (t-s) } \, \bigl\{ e^{-\beta (x-s)} |\gamma_\eps(x-s)|  \bigr\} \bigl\{ e^{-\beta(t-s)} |N(t-s) | \bigr\}\,  \mathrm{d}s
\\
&&\quad\le C \, e^{{a^*} t} \| (|\gamma_\eps| e_{-\beta} ) *  (|N |  e_{-\beta} ) \|_{C_b} 
\\
&&\quad\le C \, e^{{a^*} t} \| \gamma_\eps e_{-\beta} \|_\XX  \, \| N  e_{-\beta} \|_{C_b}.
       \end{eqnarray*}
Using the definition \eqref{def:Aeps} and the estimates \eqref{eq:F&FprimBd} and  \eqref{eq:N<ea}, we deduce 
\beqn\label{eq:SBASBLunif&expo} 
\|  (S_{\BB}\ast\AA S_{\BB})(t)g \,  e_{-\beta} \|_{L^\infty} 
\le C' \,e^{{a^*} t} \,  \| g \|_\XX, \quad \forall \, t \ge 0, 
\eeqn
for some constant $C' \in (0,\infty)$. 

\smallskip
 
Next, differentiating the  functions in both sides of identity \eqref{eq:identitySBASB2}, we get 
\begin{eqnarray*}
\partial_x[ (S_{\BB}\ast\AA S_{\BB})(t)g] (x) 
&=& -a(x) e^{-A(x)}(\nu_\eps\ast\check N_t)(x)
-e^{-A(x)}(\nu_\eps\ast\check N'_t)(x)\\
   && - e^{-A(x)} \, \nu_\eps(x-t) N(0) \mathbf{1}_{x-t\geq0} +e^{-A(x)} \, \nu_\eps(x) N(t),
  \end{eqnarray*}
  with $\check N'_t(s)=N'(t-s)$. In order to estimate the second term, we compute
 \bean
  N'(t)
  &=& (1-\kappa)^{-1}\int_0^\infty\partial_x[a(x)e^{-A(x)}] \, e^{A(x-t)} g(x-t)\mathbf{1}_{x-t\geq0}\,\mathrm{d}x
  \\
  &=& (1-\kappa)^{-1}\int_0^\infty [a'(x) - a(x)^2] \,  e^{A(x-t)-A(x)} g(x-t)\mathbf{1}_{x-t\geq0}\,\mathrm{d}x.
  \eean
Using the inequality \eqref{eq:eA-A}, we deduce
\bean
|N'(t)|
&\le& C' \int_0^\infty\!\!\bigl\{ |a'(x)| + a(x)^2 \bigr\}\, e^{ 3\beta t}   g(x-t)\mathbf{1}_{x-t\geq0}\,\mathrm{d}x dt
\\
&\lesssim& e^{{a^*} t}  \| g \|_\XX.
\eean 
As a consequence, using again  \eqref{eq:eA-A}, we get 
\bean
  \| e^{-A }(\nu_\eps\ast\check N'_t) \|_\XX 
& \lesssim&
  \int_0^t \! \!\int_0^\infty e^{3\beta s} \, |\gamma_\eps(x-s)| \, |N'(t-s)| \, dxds 
\\
& \lesssim&
 e^{{a^*} t} \,   \| \gamma_\eps \|_\XX \, \| N'   \|_\XX \lesssim   e^{{a^*} t} \, \| g \|_\XX. 
\eean
Treating  in a similar way the two other terms involved in the expression of $\partial_x[ (S_{\BB}\ast\AA S_{\BB})(t)g]   $, we finally obtain 
\bean
  \| \partial_x[ (S_{\BB}\ast\AA S_{\BB})(t)g]  \|_\XX  \lesssim e^{{a^*} t} \, \| g \|_\XX. 
\eean
Together with \eqref{eq:SBASBLunif&expo} that concludes the proof of (iii). 
\Black
\endproof

\proof[Proof of Theorem~\ref{th:MRe}] The main idea is to apply or adapt the versions of Weyl's Theorem  \cite[Theorem 3.1]{MS},
\cite[Theorem E.3.1]{MS-erratum} and of the Spectral Mapping Theorem \cite[Theorem 2.1]{MS}, \cite[Theorem E.2.1]{MS-erratum} (see also the variant results in \cite{Mbook*}). 
We start collecting the three key properties satisfied by the operators $\AA_\eps$ and $\BB_\eps$ involved in the splitting of $\Lambda_\eps$ 

\smallskip
We denote by $X_\zeta$ the abstract Sobolev space associated to $L_\eps$ for $\zeta \in \R$, see e.g. \cite[Section~II.5]{EngelNagel}, so that in particular $X_\zeta \subset W^{\zeta,1}(\R_+)$, with equality when $\zeta \le 0$. We recall that for the generator $L$ of a semigroup $S_L$, we define the resolvent set $\rho(L)$   by 
$$
\rho(L) := \{ z \in \C; \,\, L-z \hbox{ is a bijection} \, \} = \C \backslash \Sigma(L), 
$$
as well as the resolvent (operator) $R_L(z) := (L-z)^{-1}$ for any $z \in \rho(L)$. We finally define the half plane $\Delta_a := \{ z \in \C; \, \Re e z > a \}$ for any $a \in \R$.

\begin{enumerate}
  \item[{\bf (A1)}] 
For any $a>{a^*}$
  and $\ell\in\mathbb{N}$, there exists  a positive constant $C_{a,\ell}$ such that the following growth estimate holds
  \begin{equation}\label{eq:A1}
    \|S_{\BB_\eps} \ast(\AA_\eps S_{\BB_\eps} )^{(\ast\ell)}(t)\|_{\mathscr{B}(X)}
    \leq C_{a,\ell} \, e^{at}, \ \ \forall t\geq0.
  \end{equation}
  It is obvious that $\mathbf{(A1)}$ is true for $\ell=0$ from Lemma~\ref{lem:Lm1}$-(ii)$ and  for $\ell = 1$ from \eqref{eq:SBASBLunif&expo}. 
  We then deduce $\mathbf{(A1)}$ for any $\ell \ge 2$ and $a > {a^*}$ by writing 
  $S_{\BB_\eps} \ast(\AA_\eps S_{\BB_\eps} )^{(\ast\ell)} = [ S_{\BB_\eps} \ast(\AA_\eps S_{\BB_\eps} ) ]  \ast(\AA_\eps S_{\BB_\eps} )^{(\ast\ell-1)}$ and by using that $(\BB_\eps-{a^*})$ is hypodissipative in $\XX$ and $\AA_\eps\in\mathscr{B}(\XX)$.
   
  \item[{\bf (A2)}]  There holds $\AA_\eps \in \BBB(X_{-1},X_\zeta)$ for any $\zeta \in (-1,0)$ and the family of operators $(S_{{\BB_\eps}} \AA_\eps  )^{(\ast 2)}(t)$ satisfies the estimate
  \begin{equation}\label{eq:A2}
  \Blue  \| (S_{\BB_\eps} \AA_\eps  )^{(\ast2)}(t)\|_{\mathscr{B}(X_{-1},X)}\leq C'_{{a^*},1}e^{{a^*} t}, \ \ \forall
    t\geq0, 
  \end{equation}
for a  positive constant $C'_{a,1}$. 
  The first claim is a consequence of  the continuous embedding $\XX \subset W^{\zeta,1}$ and of Lemma~\ref{lem:Lm1}$-(i)$. The
  second estimate is obtained by putting together the properties $(i)$ and $(ii)$ in Lemma~\ref{lem:Lm1}.
 
  \item[{\bf (A3)}] \Blue The family of operators $(R_{{\BB_\eps}} \AA_\eps  )^3(z)$ satisfies the compactness estimate
$$
  \| (R_{{\BB_\eps}} (z) \AA_\eps  )^3  \|_{\BBB(\XX,Y)}  \leq C''_{2,a} , \quad \forall \, z \in \Delta_a, \,\, \forall \, a > {a^*}, 
$$
  for some positive constant $C''_{2,a}$. Observing that  $z \mapsto - (R_{{\BB_\eps}} (z) \AA_\eps  )^3$ is nothing but the  Laplace transform of the function $t \mapsto (S_{\BB_\eps} \AA_\eps  )^{(\ast3)}(t)$, that   clearly holds true thanks to Lemma~\ref{lem:Lm1}. 
\end{enumerate}
\Blue We now briefly explain how the proof goes on, referring to \cite{MS,MS-erratum,Mbook*} for more details and developments. In order to shorten notation we skip again the $\eps$ dependency.

\smallskip\noindent
{\sl Step 1. }
On the one hand, because of \eqref{eq:A2} and the fact that $R_{\BB}$ is nothing but the opposite of the Laplace transform of $S_{\BB}$, we have 
\beqn\label{eq:RBX-1X}
\|  (R_{\BB} (z) {\AA})^2 \|_{X_{-1} \to X} \le C, \quad \forall \, z \in \Delta_a.
\eeqn

Recalling that the negative abstract Sobolev norm is defined by 
$$
\| g \|_{X_{-1}}  := \| R_{\BB}(a) g \|_X,
$$
and using the dissipativity property of ${\BB}$, we immediately have 
$$
\forall \, z \in \Delta_a, \,\, \forall \, g \in X, \quad \| R_{\BB}(z) g \|_{X_{-1}} \le C \, \| g \|_{X_{-1}}.
$$
\Blue Moreover, using the resolvent identity 
\bean
R_{\BB}(a) R_{\BB}(z) 
&=& z^{-1} \, ( R_{\BB}(a)  {\BB} R_{\BB}(z) - R_{\BB}(a)),
\\
&=& z^{-1} \, ( R_{\BB}(z) - R_{\BB}(a) + a R_{\BB}(z) R_{\BB}(a)),
\eean
we also have 
\bean
\| R_{\BB}(z) g \|_{X_{-1}} 
&\le& {1 \over |z |} \, \{ \| R_{\BB}(z) g \|_X + \| R_{\BB}(a) g \|_X+ |a| \, \| R_{\BB}(z) R_{\BB}(a) g \|_X\}
\\
&\le& {C \over |z |} \, \| g \|_X ,
\eean
for any $g\in X$ and $z \in \Delta_a$.  Finally, using both estimates and an interpolation argument, we deduce 
\beqn\label{eq:RBXsX-1}
\forall \, z \in \Delta_a, \quad \| R_{\BB}(z) \|_{X_{\zeta}\to X_{-1} } \le {C \over \langle z \rangle^{1+\zeta}},  
\quad \forall \, \zeta \in (-1,0).
\eeqn
We conclude with 
\beqn\label{eq:RBXX}
\forall \, z \in \Delta_a, \quad \| (R_{\BB} (z) {\AA})^{3} \|_{\XX \to X } \le {C \over \langle z \rangle^{3/4}}, 
\eeqn
by writing 
$$
(R_{\BB} (z) {\AA})^{3} = (R_{\BB} (z) {\AA})^2 R_{\BB} (z) {\AA}
$$
and by using estimates (i) in Lemma~\ref{lem:Lm1} with $\zeta=-1/4$, \eqref{eq:RBXsX-1} with $\zeta = -1/4$ and \eqref{eq:RBX-1X}.

\smallskip\noindent
{\it Step 2. } Here, we follow \cite[Proof of Theorem~2.1]{MS} and \cite{Voigt80} by taking advantage of the analysis of 
degenerate-meromorphic functions performed in \cite{RibaricVidav}. 
  Iterating the factorization identity 
$R_{\Lambda} = R_{\BB} - R_{\BB} \AA R_{\Lambda}$, we obtain 
\beqn\label{RI-VLamb=U}
 (I - \VV)  R_{\Lambda}= \UU \quad \hbox{in} \quad \BBB(\XX)
\eeqn
with 
$$
\UU = R_{\BB} - R_{\BB} \AA R_{\BB}+ (R_{\BB} \AA)^2 R_{\BB}, \quad \VV := - (R_{\BB} \AA)^3 .
$$
Because of  \eqref{eq:RBXX}, \eqref{eq:A1} and of the identity
$$
R_{\BB} (\AA R_{\BB})^\ell (z)  = (-1)^{\ell+1}\int_0^\infty e^{-zt}  S_{\BB} \ast(\AA S_{\BB} )^{(\ast\ell)}(t) \, dt  \in \BBB(X)
 $$
for any $z \in \Delta_a$, we see that all the terms involved in  \eqref{RI-VLamb=U} act in $\BBB(X)$, so that we also have 
\beqn\label{RI-VL=U}
(I - \VV) R_{L}  = \UU \quad \hbox{in} \quad \BBB(X). 
\eeqn
Because of \eqref{eq:RBXX}, we see that $I-\VV$ is inversible on $\Delta_a \backslash B(0,M)$ for $M > 0$ large enough and  \eqref{RI-VL=U} tells us that $R_{L} = (I-\VV)^{-1} \UU$ is uniformly bounded in $\Delta_a \backslash B(0,M)$, in particular $\Sigma(L) \subset \Delta_a^c \cup B(0,M)$. 

\smallskip
On the other hand, because of  {\bf (A3)}, 
there holds 
$$
\|  \VV(z) \|_{X \to Y} \le C,
$$
with compact embedding $Y \subset X$. As a consequence of the fact that $\VV$ is the opposite of the Laplace transform of the nice function $\R_+ \to X$, $t \mapsto (S_\BB \AA)^{(*3)}(t)$, we also deduce that  $\VV$ is holomorphic in $\BBB(X)$ on $\Delta_a$. 
We define   $\Psi(z) := I - \VV(z)$. Because $\Psi$ is an holomorphic function on $\Delta_a$, $\Psi(z)$ is invertible for (some) $z \in \Delta_a \backslash B(0,M)$ and $R \VV(z) \subset Y$ for any $z \in \Delta_a$, \cite[Corollary~II]{RibaricVidav} asserts that the function $z \mapsto \Psi(z)^{-1}$ is a degenerate-meromorphic operators valued function, that is  $\Psi$ is holomorphic on $\Delta_a \backslash D$, $D \subset \Delta_a$ is discrete, any point $\xi \in D$ is an  isolate pole and the coefficients of the principal part in the Laurent series are finite rank operators. Next,   from  \cite{RibaricVidav}  and the identity  
$$
 R_{L} =  \Psi^{-1}\UU,
$$
 we deduce that $ R_{L}$ is also degenerate-meromorphic on $\Delta_a$, in particular $\Sigma(L) \cap \Delta_a$ is a finite set of discrete eigenvalues,
  which means isolated eigenvalues associated to an algebraically finite eigenspace. That is nothing but \eqref{eq:SigmaLambdaE1}. We define $\Pi := I-\Pi_{\Lambda,a}$.

\smallskip\noindent
{\sl Step 3. } We claim that for any $n \ge 2$, there holds  
\bear
\label{eq:PiSL=}
S_L(t) \Pi  &=& \{ S_\BB (t) + ... + (S_\BB  \AA )^{*(n-1)} * S_\BB(t) \}  \Pi 
\\ \nonumber
&&+ (-1)^{n+1} {i \over 2\pi} \int_{\uparrow_a} (R_\BB\AA )^n (z) \, R_L (z) \Pi  \, e^{t z} \, dz
\eear
in $\BB(X)$ where $\uparrow_c := \{ c + i y, \, y \in \R\}$ is the complex line of abscissa $c \in \R$. In order to prove that identity, we iterate the Duhamel formula 
$$
S_\Lambda(t)  = S_\BB(t) + S_\BB \AA * S_\Lambda (t) ,
$$
in $\BBB(\XX)$ and for any $f \in X$ and $n\ge2$. We then obtain  
$$
S_L(t) \Pi f =  \{ S_\BB (t) + ... + (S_\BB  \AA )^{*(n-1)} * S_\BB(t) \}  \Pi f + g_n(t), 
$$
with 
$$
g_n(t) := (S_\BB  \AA )^{*(n)} * S_L(t) \Pi f.
$$
Because $t \mapsto g_n'(t) e^{-bt} = \BB g_n(t) \, e^{-bt} + \AA g_{n-1}(t)\, e^{-bt}\in L^1(\R_+,\XX)$ for $b > 0$ large enough from {\bf (A1)}-{\bf (A2)}  and the fact that $t \mapsto S_L(t) \Pi f \, e^{-bt} \in L^\infty(X)$, the following inverse Laplace transform formula holds 
$$
 g_n(t) =  (-1)^{n+1} {i \over 2\pi} \int_{\uparrow_b} (R_\BB\AA )^n (z) \, R_L (z) \Pi  f \, e^{t z} \, dz, \quad \forall \, t \ge 0.
$$
From the preceding steps, the function under the integral sign is holomorphic on $\Delta_{a^*}$ and we may move the line of integration from $\uparrow_b$ to $\uparrow_a$.
Observing that every term in the resulting identity belongs to $X$, we conclude to  \eqref{eq:PiSL=}. More precisely and importantly, taking $n=6$, using \eqref{eq:RBXX} and the fact that $\| R_L (z) \Pi \|_{\BBB(X)}$ is uniformly bounded on $\Delta_a$ from Step 2, we have 
$$
\| g_n(t) \|_X \le  {e^{at} \over 2\pi} \int_{\uparrow_a}  \| (R_\BB\AA )^6 (z) \|_{\BBB(X)} \| R_L (z) \Pi \|_{\BBB(X)} \, dz \| f \|_X \lesssim e^{at} \| f \|_X,
$$
for any $t \ge 0$. We conclude to \eqref{eq:SLambdaepsPiPerp}  using that the other terms involved in  \eqref{eq:PiSL=} are similarly bounded thanks to {\bf (A1)} and  using a weakly $*$ density argument.
\Black
\qed

\subsection{The vanishing connectivity regime} 

When the network connectivity parameter vanishes, $\eps=0$, the linearized time elapsed operator simplifies into 
\begin{equation}\label{TLC}
    \Lambda_0g=-\partial_x g-a(x,0)g+\delta_{x=0}\MM_0[g],
\end{equation}
where $\MM_0[g]=\int_0^\infty a(x,0)g(x)\mathrm{d}x$.  The associated semigroup is then positive and the 
following version of the Krein-Rutman theorem holds.

\begin{theorem}\label{lc1}
There exist some constants $\alpha_0<0$ and $C>0$ such that $\Sigma(\Lambda_0) \cap \Delta_{\alpha_0} = \{ 0 \}$ and for any $g_0 \in \XX$, $\langle g_0 \rangle = 0$, there holds 
\beqn\label{eq:lc1}
\| S_{\Lambda_0}(t) g_0 \|_\XX \leq C e^{\alpha_0 t} \, \| g_0 \|_\XX, \quad \forall \, t \ge 0.
\eeqn
\end{theorem}

We denote 
$$
\XX_+ := \{ g \in M^1(\mathbb{R}_+); \,\, g \ge 0 \}, 
$$
the space  of  bounded and nonnegative Radon measures. 

\smallskip
We start with two elementary auxiliary results. 

\begin{lemma}\label{lem:S0p}
$S_{\Lambda_0}$ is positive: $S_{\Lambda_0}(t)g\in
\XX_+$ for any $g\in \XX_+$ and any $t\geq0$.
\end{lemma}

\proof  We introduce a dual problem of \eqref{TLC} defined on the space $C_0(\R)$ by  
\beqn\label{eq:tildeLambda}
\partial_t \varphi  = \tilde \Lambda \varphi    = \tilde \BB  \varphi  +  \tilde \AA \varphi ,
\quad \varphi(0,\cdot) = \varphi_0, 
\eeqn
where the operators $\tilde \AA$ and $\tilde \BB$ are defined by 
$$
\tilde\BB\varphi =  \partial_x \varphi  - a(x,0) \varphi, 
\quad 
\tilde \AA \varphi = a(x,0)\varphi(0). 
$$
A solution $\varphi$ to \eqref{eq:tildeLambda} then satisfies 
$$
\varphi(t) = S_{\tilde \BB} (t)  \varphi_0 + (S_{\tilde \BB} * \tilde\AA \varphi)(t). 
$$
Let us  fix $\varphi_0 \in C_0(\R)$ such that $\varphi_0 \le 0$ and let us prove that $\varphi(t) \le 0$ for any $t \ge 0$. 
We obviously have that $S_{\tilde \BB}$ is a positive operator and it is a contraction  in $C_0(\R)$. Taking the positive part
in  \eqref{eq:tildeLambda}, we get 
\bean
\varphi_+(t) 
&\le&  S_{\tilde \BB} (t)  \varphi_{0+} + (S_{\tilde \BB} * \tilde\AA \varphi_+)(t)
\\
&\le&   a_1 \int_0^t S_{\tilde \BB} (t-s)  \varphi_+ (0) \, \mathrm{d}s ,
\eean
so that 
$$
\| \varphi_+ (t)\|_{L^\infty} \le C \int_0^t  \| \varphi_+ (s)\|_{L^\infty} \, \mathrm{d}s.
$$
From Gr\"onwall's lemma, we deduce that $\varphi_+(t) = 0$ for any $t \ge 0$ and then $\varphi \le 0$.  We conclude by observing that $S_{\Lambda_0}$
is the dual of $S_{\tilde\Lambda}$.
\endproof
\begin{lemma}\label{lem:Mp}
$-\Lambda_0$ satisfies the following version of the strong maximum principle: 
 for any given $g\in \XX_+$ and $\mu\in\mathbb{R}$, there holds
$$
g\in D(\Lambda_0)\setminus\{0\} \ and \ (-\Lambda_0+\mu)g\geq0 \
imply \ g>0.
$$
\end{lemma}
\proof  Suppose that there holds $(-\Lambda_0+\mu)g\geq0$ for
$g$ satisfying the above conditions. It is only necessary to prove that
$g$ does not vanish in $\mathbb{R}_+$. Since $g\not\equiv0$, there
exists $x^\ast\in\mathbb{R}_+$ such that $g(x^\ast)>0$. Rewriting
the assumption as
$$
\partial_x g+\big(a(x,0)+\mu\big)g\geq\delta_{x=0}\int_0^\infty
a(x,0)g\mathrm{d}x,
$$
 we observe that
\begin{equation}\label{SMP}
\partial_x(e^{A(x,0)+\mu x}g)=e^{A(x,0)+\mu
x}\Big(\partial_x g+\big(a(x,0)+\mu\big)g\Big) \ge 0.
\end{equation}
\begin{enumerate}
 \item[(i)] For $x\in(x^\ast,\infty)$, we have
$$
 e^{A(x,0)+\mu x}g \geq e^{A(x^\ast,0)+\mu x^\ast}g(x^\ast)>0.
$$
 \item[(ii)] For $x\in(0,x^\ast)$, by integrating the same equation on
 $(0,x)$, we obtain
 \begin{eqnarray*}
 e^{A(x,0)+\mu x}g &\geq& \int_0^x\delta_{u=0}e^{A(u,0)+\mu
 u}\int_0^\infty a(y,0)g(y)\mathrm{d}y\mathrm{d}u+g(0) \\
 &\geq& \int_0^\infty a(y,0)g(y)\mathrm{d}y.
 \end{eqnarray*}
 From the positivity assumption \eqref{hyp:a2} on $a$ and step  $\mathrm{(i)}$, we have
 $$
 \int_0^\infty
 a(y,0)g(y)\mathrm{d}y>\frac{a_0}{2}\int_{\max\{x_0,x^\ast\}}^\infty g(y)\mathrm{d}y>0.
 $$
\end{enumerate}
Therefore, $g$ does not vanish on $(0,\infty)$.
\endproof

\proof[Proof of Theorem~\ref{lc1}] First, we know from Theorem~\ref{th:SS} that there exists at least one nonnegative and non-vanishing solution $F_0$ to the
eigenvalue problem $\Lambda_0 F_0 = 0$ and the associated dual eigenvector is $\psi = 1$. 
Next, we observe that, defining the $\hbox{\rm sign} f$ operator for $f \in D(\Lambda^2_0)$ by 
$$
[(\hbox{\rm sign} f )^* \, \psi ] (x) := {1 \over 2 |f(x)|} \,  [\bar f(x) \psi(x) +  f(x) \bar \psi(x) ], \quad \forall \, \psi \in C_0(\R),
$$
we have,  for any $\psi \in C_0(\R)_+$, 
\bean
\Re e \langle(\hbox{\rm sign} f )  \MM_0[f]\, ,\psi\rangle 
&=&  \Re e \langle\MM_0[f]\, , (\mathrm{sign}f)^\ast\psi\rangle
\\
&=& \Re e \,\Bigl[ \int a_0 f \, dx \Bigr] \frac{ \Re e f(0)}{|f(0)|} \, \psi(0) 
\\
 &\leq& \int a_0|f|dx \,  \psi(0)=\langle\MM_0[|f|]\, ,\psi\rangle,
\eean
%
%
%
%
which is nothing but the complex Kato's inequality
\beqn\label{eq:KatoIneq}
\forall \, f \in D(\Lambda^2_0), \quad \Re e  (\hbox{\rm sign} f ) \, \Lambda_0 f \le \Lambda_0 |f| .
\eeqn
We also observe that $D(\Lambda_0^2) \subset C_b(0,\infty)$, and, as a consequence, $g \in D(\Lambda_0^2)$ and $|g| > 0$ implies $g > 0$ or $g < 0$. 
We then may use exactly the same argument as in \cite[Proof of Theorem~5.3]{MS,MS-erratum} (see also \cite{Mbook*}):

- Kato's inequality \eqref{eq:KatoIneq} and the strong maximum principle imply that the eigenvalue $\lambda = 0$ is simple and the associated eigenspace is $\hbox{Vect}(F_0)$; 

- together with the fact that $S_{\Lambda_0}$ is a positive semigroup, one deduces that $\lambda = 0$ is the only eigenvalue with nonnegative real part. 

\smallskip
We conclude to the spectral gap estimate \eqref{eq:lc1} for some $\alpha_0 \in ({a^*}, 0)$ with the help of Theorem~\ref{th:MRe}. 
\endproof

\subsection{Weak connectivity regime - exponential stability of the linearized equation}
We extend the exponential stability property which holds for a vanishing connectivity to the weak connectivity regime thanks to a perturbation argument. 
 
 \begin{theorem}\label{th:CWR}
There exist some constants $\eps_0 > 0$, $\alpha<0$ and $C>0$ such that for any $\eps \in [0,\eps_0]$ there hold $\Sigma(\Lambda_\eps) \cap \Delta_{\alpha} = \{ 0 \}$ and  
\beqn\label{eq:lc}
\| S_{\Lambda_\eps}(t) g_0 \|_\XX \leq C e^{{\alpha} t} \, \| g_0 \|_\XX, \quad \forall \, t \ge 0,
\eeqn
for any $g_0 \in \XX$, $\langle g_0 \rangle = 0$. 
\end{theorem}

The proof follows the stability theory for semigroups developed in Kato's book \cite{Kato} and revisited in \cite{MMcmp,Tristani,Mbook*,MT*}. 

\Blue
\proof[Proof of Theorem~\ref{th:CWR}] 
With the definitions \eqref{def:Nepsg},  \eqref{def:Aeps} and \eqref{def:Beps} of $\MM_\eps$, $\AA_\eps$  and $\BB_\eps$, we have 
$$
(\BB_\eps - \BB_0 )g = (a(x,0)-a(x,\eps M_\eps ))g 
$$
and 
$$
(\AA_\eps - \AA_0 )g = (\MM_\eps[g] - \MM_0[g]) \, \delta_0 -  \eps (\partial_\mu a)(x,\eps M_\eps)  \, F_\eps \, \MM_\eps[g].
$$
Together with the smoothness assumption \eqref{hyp:a3}, we deduce
\beqn\label{eq:Lambdaeps-0}
\| \BB_\eps - \BB_0 \|_{\BBB(\XX)} + \| \AA_\eps - \AA_0 \|_{\BBB(\XX)} \le C \, \eps, \quad \forall \, \eps \ge 0.
\eeqn
We then argue similarly as in the proof of \cite[Theorem 2.15]{Tristani} (see also \cite{Kato,MMcmp,Mbook*,MT*}) and therefore  just sketch the proof. 
We now define 
$$
K_\eps(z) :=    (R_{\BB_\eps} (z) \AA_\eps)^2 R_{\Lambda_0}(z)  (\Lambda_\eps - \Lambda_0) \in \BBB(\XX,X)
$$
and we take any $\alpha \in (\alpha_0,0)$, recalling that $a^* < \alpha_0 < 0$. We deduce from \eqref{eq:Lambdaeps-0} and the estimates (i) and (ii) in Lemma~\ref{lem:Lm1} that for some $\eta < |\alpha|$, $\eps_0 > 0$, small enough, $C >0$, and for any $z \in \Delta_\alpha \backslash B(0,\eta)$, $\eps \in [0,\eps_0)$,
we have 
\beqn\label{eq:Keps-0}
\| K_\eps (z) \|_{\BBB(X)} \le C \eps < C \eps_0 < 1, 
\eeqn
and thus $(1-K_\eps(z))^{-1}$ is well defined in $\BBB(X)$. Using the elementary identities
\beqn\label{eq:RL=RB2}
R_{\Lambda_\eps} = R_{\BB_\eps} - R_{\BB_\eps} \AA_\eps R_{\BB_\eps} +  (R_{\BB_\eps} \AA_\eps)^2 R_{\Lambda_\eps}  
=: \UU_\eps + (R_{\BB_\eps} \AA_\eps)^2 R_{\Lambda_\eps} , 
\eeqn
and 
$$
R_{\Lambda_\eps} = R_{\Lambda_0} + R_{\Lambda_0} (\Lambda_\eps - \Lambda_0) R_{\Lambda_\eps}, 
$$
we get
$$
 (I - K_\eps) R_{\Lambda_\eps} = \UU_\eps  + (R_{\BB_\eps} \AA_\eps )^2 R_{\Lambda_0} .
$$
Because both terms are well defined in $\BBB(X)$ from Lemma~\ref{lem:Lm1} and the link between resolvent and semigroup, we also 
have 
$$
 (I - K_\eps) R_{L_\eps} = \UU_\eps  + (R_{\BB_\eps} \AA_\eps )^2 R_{L_0},
$$
from what we deduce
$$
R_{L_\eps} =  (I -  K_\eps)^{-1}  ( \UU_\eps  + (R_{\BB_\eps} \AA_\eps )^2 R_{L_0}  ).
$$
Since the RHS expression is clearly uniformly bounded in $\BBB(X)$ on the complex set $\Delta_\alpha \backslash B(0,\eta)$ for any $\eps \in [0,\eps_0)$, we have $\Sigma(L_\eps) \cap \Delta_\alpha \subset B(0,\eta)$. 

\smallskip
Using the definition of the eigenprojector $\Pi_0$ on the eigenspace associated to the spectral values of $\Lambda_0$ lying in $B(0,\eta)$ by mean of Dunford integral (see \cite[Section III.6.4]{Kato} or \cite{GMM*,Mbook*}) and the analogous of \eqref{eq:RL=RB2} in $\BBB(X)$, we have 
\bean
\Pi_0&=& {i \over 2\pi} \int_{|z| = \eta} R_{L_0}(z) \, dz
\\
&=& {i \over 2\pi} \int_{|z| = \eta}   (R_{\BB_0} \AA_0 )^2 R_{L_0} \, dz, 
\eean
by using that the contribution of holomorphic functions vanish.
In a similar way, we have 
\bean
\Pi_\eps 
&=&  {i \over 2\pi} \int_{|z| = \eta} (I-K_\eps + K_\eps) R_{L_\eps} \, \mathrm{d}z
\\
&=& {i \over 2\pi} \int_{|z| = \eta}   (R_{\BB_\eps} \AA_\eps )^2 R_{L_0} \, dz
+  {i \over 2\pi} \int_{|z| = \eta}  K_\eps R_{L_\eps}  \, dz.
\eean
For $g \in X$, we next compute
\bean
\| (\Pi_\eps - \Pi_0) g \|_{X} 
&\le& {1 \over 2\pi}  \int_{|z| = \eta}   \| ((R_{\BB_\eps} \AA_\eps )^2  -  (R_{\BB_0} \AA_0 )^2) R_{\Lambda_0} g \|_\XX \, dz
\\
&& +  {1 \over 2\pi}  \int_{|z| = \eta}   \| K_\eps R_{\Lambda_\eps} g \|_X \, dz \le C \, \eps \, \| g \|_\XX,
\eean
where we have used the  identity 
\bean
(R_{\BB_\eps} \AA_\eps )^2  -  (R_{\BB_0} \AA_0 )^2
&= &R_{\BB_\eps} \AA_\eps R_{\BB_\eps}  \{ (\AA_\eps- \AA_0 ) + (\BB_0 - \BB_\eps) R_{\BB_0} \AA_0\}
\\
&&+ R_{\BB_\eps}  \{ (\AA_\eps- \AA_0 ) +  (\BB_0 - \BB_\eps) R_{\BB_0} \AA_0 \} R_{\BB_0} \AA_0 
\eean
and the estimates \eqref{eq:Keps-0} and \eqref{eq:Lambdaeps-0}. As a consequence, we deduce 
$$
\| \Pi_\eps - \Pi_0 \|_{\BBB(X)} < 1.
$$
for any $\eps \in (0,\eps_0)$, up to take a smaller real number $\eps_0 > 0$.
From the classical result \cite[Section I.4.6]{Kato} (or more explicitly \cite[Lemma~2.18]{Tristani}), we deduce that there exists $\xi_\eps \in \Delta_\alpha$ such that 
$$
\Sigma(\Lambda_\eps)\cap\Delta_\alpha=\{\xi_\eps\}, \quad \xi_\eps \hbox{ is a simple eigenvalue}, 
$$
for any $\eps \in [0,\eps_0]$. We conclude by observing that $\xi_\eps = 0$ because $1 \in X'$ and $\Lambda_\eps^* 1 = 0$ (which is nothing but the mass conservation). 
\qed
\Black


%
\subsection{Weak connectivity regime - nonlinear exponential stability}

Now, we focus on the nonlinear exponential stability of the solution to the evolution equation   \eqref{ASM}--\eqref{eq:ASM3} in the case without delay.
We start with an auxiliary result. 
We define the function $\Phi :  L^1(\R_+) \times \R \to \R$ by 
$$
\Phi [g,\mu] := \int_0^\infty a(x,\eps \mu) g (x) \, \mathrm{d}x - \mu.
$$
We denote by $W_1$ the optimal transportation Monge-Kantorovich-Wasserstein distance on the probability measures set $\Ppp(\R_+)$ associated to the distance $d(x,y)= |x-y| \wedge 1$, or equivalently defined by 
$$
\forall \, f,g\in \Ppp(\R_+), \quad W_1(f,g) := \sup_{\varphi, \| \varphi \|_{W^{1,\infty}} \le 1} \int_0^\infty (f-g) \, \varphi.
$$


\begin{lemma}\label{lem:varphig} Assume \eqref{hyp:a3}. There exists $\eps_0 > 0$ and for any $\eps \in (0,\eps_0)$ there exists a function $\varphi_\eps : \Ppp(\R) \to \R$ which is Lipschitz continuous for the weak topology of probability measures and such that $\mu = \varphi_\eps[g]$ is the unique solution to the equation  
$$
 \mu \in \R_+, \quad \Phi(g,\mu) = 0.
$$
\end{lemma}

\noindent{\sl Proof of Lemma~\ref{lem:varphig}. } {\sl Step 1. Existence. }
For any $g \in \Ppp(\R)$ we have $\Phi(g,0) > 0$, and for any $g \in \Ppp(\R)$ and $\mu \ge 0$, we have 
$$
 \Phi(g,\mu) \le \| a \|_{L^\infty}  -  \mu,
$$
so that $\Phi(g,\mu)<0$ for $\mu>\|a\|_{L^\infty}$. By the intermediate value theorem and the continuity property of $\Phi$, for any fixed $g \in \Ppp(\R_+)$ and $\eps \ge 0$, there exists at least one solution $\mu \in (0,\| a \|_{L^\infty}]$ to the  equation $ \Phi(g,\mu) = 0$. 
 
\smallskip\noindent {\sl Step 2. Uniqueness and Lipschitz continuity. }
 Fix $f,g \in \Ppp(\R_+)$ and consider $\mu,\nu \in \R_+$ such that 
$$
 \Phi(f,\mu) =  \Phi(g,\nu) = 0.
 $$
We have  
$$
\nu - \mu = \int_0^\infty a(x,\eps \nu) (g - f) +   \int_0^\infty (a(x,\eps \nu) -  a(x,\eps \mu) ) f ,
$$
with 
$$
\Bigl| \int_0^\infty a(x,\eps \nu) (g - f) \Bigr| \le \| a (\cdot,\eps\nu) \|_{W^{1,\infty}} \, W_1(g,f),
$$
and
$$
\Bigl| \int_0^\infty \big(a(x,\eps \nu) -  a(x,\eps \mu)\big) f \Bigr| \le \| a(\cdot ,\eps \nu) -  a(\cdot,\eps \mu) \|_{L^{\infty}} 
\le \eps\, \| \partial_\mu a \|_{L^\infty} |\mu-\nu|.
$$
We then obtain 
\beqn\label{eq:mu-nu}
|\mu - \nu| \, (1-\eps  \| \partial_\mu a \|_{L^\infty}) \le \| a (\cdot,\eps\nu) \|_{W^{1,\infty}} \, W_1(g,f),
\eeqn
and we may fix $\eps_0 > 0$ such that $1-\eps_0  \| \partial_\mu a \|_{L^\infty} \in (0,1)$, $\eps \in [0,\eps_0]$. On the one hand, for $f=g$, we deduce 
that $\mu=\nu$ and that uniquely defines the mapping $\varphi_\eps[g] := \mu$. On the other hand, the function is Lipschitz continuous because of \eqref{eq:mu-nu}.
\qed

%
%

\smallskip
We also recall the following classical Gr\"onwall's type lemma. 

\begin{lemma}\label{lem:Gronwall} Assume  that $u \in C([0,\infty); \R_+)$ satisfies the integral inequality 
$$
u(t) \le C_1 e^{a t} \, u_0 + C_2 \int_0^t e^{a (t-s)} u(s)^2 \, \mathrm{d}s, \quad \forall \, t > 0, 
$$
for some constants $C_1 \ge 1$, $C_2, u_0 \ge 0$ and $a < 0$. Under the smallness assumption 
$$
a + 2 C_2 u_0 < 0, 
$$
there holds 
$$
u(t) \le \Bigl( 1 + {C_1 u_0 C_2 \over |a + 2   C_2 u_0|} \Bigr) C_1 \, e^{at} \, u_0, \quad \forall \, t \ge 0.
$$

\end{lemma}

\noindent{\sl Proof of Lemma~\ref{lem:Gronwall}. } We fix $A \in ( C_1 u_0,2C_1u_0)$, so that $C_1 u(t) \le A$ at least on a small interval, that is for any $t \in [0,\tau]$, $\tau > 0$ small enough, 
and then  the integral inequality implies on the same interval 
$$
u(t) \le C_1 e^{a t} \, u_0 + C_2 C_1^{-1} A \, \int_0^t e^{a (t-s)} u(s) \, \mathrm{d}s.
$$
The classical  Gr\"onwall's lemma (for linear  integral inequality) and the smallness assumption $a + C_2 C_1^{-1} A \le 0$ imply  
$$
u(t) \le C_1 \, u_0 \, e^{(a+C_2  C_1^{-1} A) t} \le  C_1 \, u_0  < A
$$
on that interval. By a continuity argument, the first above  inequality holds on $\R_+$ and then with $A := C_1 u_0$.  Next, replacing that first estimate in the integral inequality we started with, we get 
$$
u(t) \le C_1 e^{a t} \, u_0 + C_2 C_1^2 u_0^2 e^{at} \int_0^t  e^{(a+2C_2  u_0) s} \, \mathrm{d}s, \quad \forall \, t > 0, 
$$
from which we immediately conclude. \qed

\medskip 
We come to the proof of our main result Theorem~\ref{th:MR} in the case without delay.
 
\proof[Proof of Theorem~\ref{th:MR} in the case without delay] We split the proof into two steps. 

\noindent
{\sl Step 1. New formulation.} 
We start giving a new formulation of the solutions to the evolution and stationary equations
in the weak connectivity regime $\eps \in (0,\eps_0]$, where $\eps_0$ is defined in Lemma~\ref{lem:varphig}. 
For a given initial datum $0 \le f_0 \in L^1(\R_+)$ with unit mass the solution $f \in C([0,\infty); L^1(\R_+))$ to the evolution
equation  \eqref{ASM} and the solution $F_\eps$ to the stationary equation \eqref{eq:StSt} clearly satisfy
\bean
\partial_t f + \partial_x f +a(\eps \varphi[f]) f=0,  & & \quad f(t,0) =  \varphi[f(t,\cdot)], 
\\
 \partial_x F +a(\eps M) F=0, & & \quad  F (0) =  M = \varphi[F],
\eean
where here and below the $\eps$ and $x$ dependency is often removed without risk of misleading.

\smallskip
We introduce the variation function $g := f - F$ which satisfies the PDE
\bear
\partial_t g 
=-  \partial_x g - a(\eps M) g - \eps a'(\eps M) F \, \MM[g]  - Q[g]
\eear
with 
$$
Q[g] :=   a(\eps \varphi[f]) f - a(\eps \varphi[F]) F  - a(\eps  \varphi[F]) g - \eps a'(\eps  \varphi[F]) F \, \MM[g],  
$$
where $\MM = \MM_\eps$ is defined in \eqref{def:Nepsg}. 
The above PDE is complemented with the boundary condition
\bean
g(t,0) &=&  \varphi[f(t,\cdot)] - \varphi[F],
\eean
and we may write again
\bean
\varphi[f]  - \varphi[F]
&=& \int_0^\infty a(\eps \varphi[f]) f  - \int_0^\infty a(\eps \varphi[F]) F 
\\
&=& \int_0^\infty \big(a(\eps M) g + \eps a'(\eps M) F \, \MM[g]\big)
+ \int_0^\infty Q[g] \, \mathrm{d}x 
\\
&=&  \MM[g] +  \QQ[g], \qquad \QQ[g] := \langle Q[g] \rangle. 
\eean
As a consequence, we have proved that the variation function $g$ satisfies  the equation 
\beqn\label{eq:dtg=Lambda+Z}
\partial_t g = \Lambda_\eps g + Z[g], \quad Z[g] :=  - Q[g] + \delta_0 \QQ[g].
\eeqn

\smallskip
\noindent{\sl Step 2. The nonlinear term. } On the one hand, we obviously have 
\beqn\label{eq:<Z>=0}
\langle Z[g ] \rangle = 0, \quad \forall \, g \in M^1(\R_+).
\eeqn
On the other hand, in order to get an estimate on the nonlinear term $Z[g]$, we introduce the notation
$$
\psi(u)=a(x, \eps m_u)f_u, 
$$
where, for some fixed $g \in \Ppp(\R_+)$, $\langle g \rangle = 0$, we  have set
$$
 f := F + g, \quad  f_u :=uf+(1-u)F, \quad m_u := \varphi[f_u].
 $$
We first notice that $\psi(0)=a(\eps \varphi[F]) F$ and $\psi(1)=a(\eps \varphi[f]) f$.
Second, we have 
\beqn\label{eq:psi'u}
\psi'(u) = a'_\eps ( m_u) f_u \, m'_u + a_\eps(m_u) g. 
\eeqn
In order to compute $m'_u$, we differentiate with respect to $u$ the identity 
$$
m_u = \int_0^\infty a_\eps(m_u) f_u  \mathrm{d}x,
$$
and we have
$$
m'_u =  \int_0^\infty a'_\eps( m_u) f_u  \mathrm{d}x \, m'_u 
+\int_0^\infty a_\eps( m_u) \, g \mathrm{d}x,
$$
which implies
\beqn\label{eq:m'u}
m'_u = \bigl( 1 -  \int_0^\infty a'_\eps( m_u) f_u  \mathrm{d}x \bigr)^{-1}   \int_0^\infty a_\eps( m_u) \, g \mathrm{d}x. 
\eeqn
We may thus observe that $m'_0 = \MM[g]$, so that  $\psi'(0)=a'_\eps( M )F \MM_\eps[g] + a_\eps( M) g$, and therefore
$$
Q[g] = \psi (1) - \psi(0) - \psi'(0).
$$
Third, from \eqref{eq:psi'u}, we have 
$$
\psi''(u) =  a''_\eps ( m_u) f_u \, (m'_u)^2  + 2  a'_\eps( m_u) g m'_u+  a'_\eps(  m_u) f_u \, m''_u, 
$$
and from  \eqref{eq:m'u}, we have 
\bean
m''(u) &=&
2 \Bigl( 1 - \int_0^\infty a'_\eps f_u \Bigr)^{-2} \, \int_0^\infty a_\eps g \int_0^\infty a'_\eps g 
\\
&&+ 2 \Bigl( 1 - \int_0^\infty a'_\eps f_u \Bigr)^{-3} \, \int_0^\infty a''_\eps f \Bigl( \int_0^\infty a_\eps g  \Bigr)^2. 
\eean
In the small connectivity regime $\eps \in (0,\eps_0]$, $\eps_0 \| a' \|_\infty < 1$, we get the bound 
\bean
\| \psi''(u) \|_X 
&\le& \| a''_\eps\|_\infty \, |m'_u|^2  + 2  \| a'_\eps \|_\infty  \| g \|_X |m'_u|  + \| a'_\eps \|_\infty \, |m''_u|
\\
&\le&  \eps^2   {\| a'' \|_\infty \| a \|_\infty^2 \over (1 - \eps \|a'\|_\infty)^2 } \, \| g \|_X^2  + 
2 \eps {  \| a' \|_\infty \| a \|_\infty \over 1 - \eps \|a'\|_\infty } \, \| g \|_X^2 
\\
&& + 2 \eps^2   {\| a' \|_\infty^2 \| a \|_\infty \over (1 - \eps \|a'\|_\infty)^2 } \, \| g \|_X^2  + 
2 \eps^3 {  \| a'' \|_\infty\| a' \|_\infty \| a \|_\infty \over (1 - \eps \|a'\|_\infty)^3 } \, \| g \|_X^2 
\\
&\le&  \eps  \, K \,  \| g \|_X^2, 
\eean
for some constant $K \in (0,\infty)$. Using  the Taylor expansion 
$$
Q[g] =\psi(1)-\psi(0)-\psi'(0)=\int_0^1(1-u)\psi''(u)\mathrm{d}u ,
$$
we then obtain 
$$
\| Z[g] \|_X \le 2 \| Q[g] \|_X \le \int_0^1 (1-u)\| \psi''(u) \|_X\mathrm{d}u \le C \, \| g \|_X^2.
$$

\smallskip
\noindent{\sl Step 3. Decay estimate. } Thanks to the Duhamel formula, the solution $g$ to the evolution equation \eqref{eq:dtg=Lambda+Z}
satisfies
$$
g(t) = S_{\Lambda_\eps}(t) (f_0 - F) + \int_0^t S_{\Lambda_\eps}(t-s)  Z[g(s)] \, \mathrm{d}s.
$$
Using Theorem~\ref{th:CWR} and the second step, we deduce  
\bean
\| g(t) \|_X
 &\le&  C \, e^{\alpha t} \, \| g_0 \|_X + \int_0^t C \, e^{\alpha (t-s) } \, \| Z[g(s)] \|_X  \, \mathrm{d}s
\\
 &\le&  C \, e^{\alpha t} \, \| g_0 \|_X + C \, \eps \, K  \int_0^t  e^{\alpha (t-s) } \, \| g(s) \|^2_X  \, \mathrm{d}s,
 \eean
for any $t \ge 0$ and for some constant $C \ge 1$, $\alpha < 0$, independent of $\eps \in (0,\eps_0]$. 
Observing that $\| g(t) \|_X \in C([0,\infty)$, we conclude thanks to Lemma~\ref{lem:Gronwall}.
\endproof


\section{Case with delay}
\label{sec:WithDelay}

This section is devoted to the proof of our main result,
Theorem~\ref{th:MR}, in the case with delay by following the same strategy as in the case without delay but adapting the functional framework.
The main difference comes from the boundary term and it will be explained in the first subsection.
We have already proved in Theorem~\ref{th:SS} the existence of a unique stationary
solution $(F_\eps,M_\eps)$ in the weak connectivity regime and we may then focus on the evolution equation.

\subsection{Linearized equation and structure of the spectrum} 
In order to write as a time autonomous equation the linearized equation \eqref{eq:ASMlin1}-\eqref{eq;ASMlin2}-\eqref{eq;ASMlin3},
we introduce the following intermediate evolution equation on a function $v=v(t,y)$
\begin{equation}\label{Eov}
    \partial_t v+\partial_y v=0, \ \
    v(t,0)=q(t), \ \ v(0,y)=0,
\end{equation}
where $y\geq0$ represents the local time for the network activity. 
That last equation can be solved with the characteristics method
$$
v(t,y)=q(t-y)\mathbf{1}_{0\leq y\leq t}.
$$
Therefore, equation \eqref{eq;ASMlin3} on the variation $n(t)$ of network activity writes
$$
n(t)= \DD[v(t)], \quad \DD[v] := \int_0^\infty v(y)b(\mathrm{d}y),
$$
and then equation \eqref{eq;ASMlin2} on the variation $q(t)$ of discharging neurons writes
$$
q(t)= \OO_\eps[g(t),v(t)] ,
$$
with 
$$
\OO_\eps[g,v] := \NN_\eps[g] + \kappa_\eps  \, \DD[v], 
$$
$$
\NN_\eps[g]  := \int_0^\infty a_\eps(M_\eps) g \, \mathrm{d}x, \quad \kappa_\eps := \int_0^\infty a'_\eps (M_\eps) F_\eps \, \mathrm{d}x.
$$
All together, we may rewrite the linear system \eqref{eq:ASMlin1}-\eqref{eq;ASMlin2}-\eqref{eq;ASMlin3}, as the autonomous system
\Blue
\begin{equation}\label{AS}
\partial_t(g,v)=\LLL_\eps(g,v),
\end{equation}
where the operator $\LLL_\eps$ is defined by
$$
\LLL_\eps
\begin{pmatrix}
g\\v
\end{pmatrix}
:=
\begin{pmatrix}
-\partial_x g-a_\eps g-a'_\eps F_\eps \DD [v]\\
-\partial_y v 
\end{pmatrix},
$$
with domain
$$
D(\LLL_\eps):=\{(g,v)\in W^{1,1}(\R_+)\times W^{1,1}(\R_+,\omega); \ g(0)=v(0)=\OO_\eps[g,v]\},
$$
where $\omega(x) : =e^{-\delta x}$ with $\delta>0$ defining in \eqref{hyp:del}.
The associated semigroup $S_{\LLL_\eps}(t)$ acts on
$$
X:=X_1\times X_2:=L^1(\R_+)\times L^1(\R_+,\omega).
$$
Considering the boundary condition as a source term, we also introduce the semigroup $S_{\varLambda_\eps}(t)$ acting on
$$
\XX:=\XX_1\times \XX_2:=M^1(\mathbb{R}_+)\times M^1(\mathbb{R}_+,\omega)
$$
with the generator $\varLambda_\eps = (\varLambda^1_\eps,\varLambda^2_\eps)$ given by 
\bean
\varLambda^1_\eps(g,v) &:=& -\partial_x g-a_\eps g-a'_\eps
F_\eps \DD [v]+\delta_{x=0}\OO_\eps[g,v],
\\
\varLambda^2_\eps(g,v) &:=& -\partial_y v +\delta_{y=0}\OO_\eps[g,v].
\eean
In a similar way as in section~\ref{subsec:WithoutD:LinearizedEq&structure}, we have $S_{\varLambda_\eps}|_X=S_{\LLL_\eps}$.

\smallskip
\Black As a first step, we establish that the semigroup {\Blue $S_{\varLambda_\eps}$ has} a nice decomposition structure with finite dimensional principal modes 
and a fast decaying remainder term. 

\begin{theorem}\label{th:MRe1} Assume \eqref{hyp:a1}-\eqref{hyp:a2}-\eqref{hyp:a3} and \eqref{hyp:del}. 
The conclusions of Theorem~\ref{th:MRe} holds true with $\Blue a^\sharp := \max\{a^*,-\delta\}<0$. 
\end{theorem}

The result is obtained as a consequence of the Spectral Mapping  and Weyl's theory developed in \cite{MS,MS-erratum,Mbook*}
and taken over in section~\ref{subsec:WithoutD:LinearizedEq&structure}. 
For that purpose, we ontroduce the convenient  splitting of the operator {\Blue $\varLambda_\eps$ on $\XX$ as  
$\varLambda_\eps =\AA_\eps+\BB_\eps$} defined by
$$
\BB_\eps (g,v) = \begin{pmatrix}
\BB^1_\eps(g,v)\\ \BB^2_\eps(g,v)
\end{pmatrix}
=
\begin{pmatrix}
-\partial_x g-a_\eps g\\
-\partial_y v 
\end{pmatrix}
$$
and
$$
\AA_\eps (g,v) =  \begin{pmatrix}
\AA^1_\eps(g,v)\\ \AA^2_\eps(g,v)
\end{pmatrix}
=
\begin{pmatrix}
-a'_\eps F_\eps\DD[v]+\delta_{x=0}\OO_\eps[g,v]\\
\delta_{y=0}\OO_\eps[g,v]
\end{pmatrix}.
$$
{\Blue It is only necessary to establish the following adequate properties of} the operators $\AA_\eps$ and $\BB_\eps$. We skip the rest of the proof and refer to the proof of Theorem~\ref{th:MRe} for more details.

\Blue
\begin{lemma}\label{lem:Lm2}
Assume that $a$ satisfies conditions \eqref{hyp:a1}--\eqref{hyp:a3} and that $b$ satisfies \eqref{hyp:del}. 
For any $\eps \ge 0$, the operators $\AA_\eps$ and
$\BB_\eps$ satisfy :  

\begin{enumerate}
  \item[(i)] $\AA_\eps\in\BBB(W^{-1,1}(\R_+)\times W^{-1,1}(\R_+,\omega),\XX)$.
  \item[(ii)] $S_{\BB_\eps}(t)$ is $a^\sharp$-hypodissipative in both $X$ and $\XX$;
  \item[(iii)] the family of operators
  $S_{\BB_\eps}\ast\AA_\eps S_{\BB_\eps}$ satisfies
  $$
  \|(S_{\BB_\eps}\ast\AA_\eps S_{\BB_\eps})(t)\|_{\BBB(\XX,Y)}\leq Ce^{\alpha
  t},\quad\forall a>a^\sharp,
  $$
  for some constant $C_a>0$ and with $Y:=Y_1\times Y_2$, where $Y_1=BV(\R_+)\cap L_1^1(R_+)$ and $Y_2=BV(\R_+,\omega)\cap L_1^1(R_+,\omega)$.
\end{enumerate}
\end{lemma}

\proof
(i) It is an immediate consequence of the fact that $\DD \in \BBB(W^{-1,1}(\R_+,\omega),\R)$ and $\NN_\eps \in \BBB(W^{-1,1}(\R_+),\R)$ because of \eqref{hyp:a3} and \eqref{hyp:del}.\Black
  
\smallskip
(ii) Since $S_{\BB^1_\eps}$ is nothing but the semigroup $S_{\BB_\eps}$ defined in \eqref{eq:defSBeps1} which is $a^*$-dissipative thanks to Lemma~\ref{lem:Lm1}-(ii), we just have to prove the dissipativity of the translation semigroup $ S_{\BB^2_\eps}$ which is defined through the explicit formula 
  $[S_{\BB^2_\eps}(t) v](y) = v(y-t)\mathbf{1}_{y-t\geq0}$.  That follows from 
  $$
  \| S_{\BB^2_\eps}(t) v \|_{X_2} \le \int_0^\infty |v(y-t)|\mathbf{1}_{y-t\geq0} \, e^{-\delta y} \, dy \le e^{-\delta t} \| v \|_{X_2}, 
  $$
  for any $v \in X_2$ and any $t \ge 0$. 

\Blue
\smallskip
(iii) For $(g,v)\in\XX$, we have 
\bear
\AA^1 S_\BB(t)(g,v) (x) &=& \gamma(x)  D(t) +\delta_{x=0}N(t), \\
\AA^2 S_\BB(t)(g,v) (y) &=& \delta_{y=0} \kappa D(t)+\delta_{y=0}N(t),
\eear
with 
\bean
\gamma(x) &:=& \kappa \delta_{x=0} -a'(x)F_\eps(x),
\\
N(t) &:=& \NN[S_{\BB^1} (t) g] = \int_0^\infty a(x)e^{A(x-t)-A(x)}g(x-t)\mathbf{1}_{x-t\geq0} \mathrm{d}x, 
\\
D(t) &:=& \DD[S_{\BB^2}(t) v] = \int_0^\infty v(y-t)\mathbf{1}_{y-t\geq0}b(\mathrm{d}y),
\eean
where here and below the $\eps$ dependency is  removed without risk
of ambiguity. We observe that
\bean
|N(t)|  &\le& C a_1e^{3\beta t} \, \| g
\|_{\XX_1},\\
|D(t)|  &\le& C e^{\delta t} \, \| v \|_{\XX_2},
\eean
for any $t\ge0$. We then compute
\bean
N'(t) &=&
\int_0^\infty \partial_x[a(x)e^{-A(x)}]e^{A(x-t)}g(x-t) \mathbf{1}_{x-t\geq0}\,
\mathrm{d}x\\
 &=& \int_0^\infty [a'(x)-a(x)^2]e^{A(x-t)-A(x)}g(x-t) \mathbf{1}_{x-t\geq0}\,
\mathrm{d}x,\\
D'(t) &=& \int_0^\infty v(y-t)  \mathbf{1}_{y-t\geq0}\, b'(\mathrm{d}y) = b' *
\check v(t),
\eean
from what we get the estimates
\bean
|N'(t)| &\lesssim& e^{a^*t}\|g\|_{\XX_1}, \\
|D'(t)| &\lesssim& e^{-\delta t}\|v\|_{\XX_2}.
\eean
Denoting 
\bean
T_1(t) (g,v) (x) &:=& (S_{\BB^1}\ast\AA^1 S_\BB)(t)(g,v)(x),\\
T_2(t) (g,v) (y) &:=& (S_{\BB^2}\ast\AA^2 S_\BB)(t)(g,v)(y), 
\eean
 we compute
\bean
   T_1(t)(g,v)(x)
   &=& \int_0^t S_{\BB^1}(s)\big(\gamma(x)D(t-s)+\delta_{x=0} N(t-s)\big)\mathrm{d}s \\
   &=& \int_0^t  e^{A(x-s)-A(x)} \big(\gamma(x-s)D(t-s)+\delta_{x-s=0}N(t-s)\mathbf{1}_{x-s\geq0}\big)\mathrm{d}s\\
   &=& e^{ -A(x)}(\nu *\check D_t)(x)+e^{-A(x)}N(t-x)\mathbf{1}_{x\le t},\\
   T_2(t)(g,v)(y)
   &=& \int_0^t S_{\BB^2}(s)\delta_{y=0}\big(\kappa  D(t-s)+N(t-s)\big)\mathrm{d}s \\
   &=& \int_0^t \delta_{y-s=0}\big(\kappa  D(t-s)+N(t-s)\big)\mathrm{d}s\\
   &=& \big(\kappa  D(t-y)+N(t-y)\big)\mathbf{1}_{y\le t},
\eean
where we use the notation $\nu  := \gamma  \, e^{A}$. 
We next differentiate the above identity, and we get 
\bean
   \partial_x T_1(t)(g,v)(x) &=& -a(x)e^{-A(x)}\big((\nu *\check
   D_t)(x)+N(t-x)\mathbf{1}_{x\le t}\big)\\
   && -e^{-A(x)}\big(\nu *\check
   D'_t(x)+N'(t-x)\mathbf{1}_{x\le t}\big)\\
   &&-e^{-A(x)}\nu (x-t)D(0)\mathbf{1}_{x-t\geq0}-e^{-A(x)}N(0)\delta_{x=t}\\
   &&+e^{-A(x)}\nu (x)D(t),\\
   \partial_y T_2(t)(g,v)(y) &=& - (\kappa  D' +N')(t-y)\mathbf{1}_{y\le t} - (\kappa  D +N)(0) \delta_{y=t}
\eean 
All together, we deduce
\bean
  \|\partial_x T_1(t)(g,v)(x)\|_{X_1} &\lesssim& e^{a^\sharp t}\|(g,v)\|_\XX, \\
  \|\partial_y T_2(t)(g,v)(y)\|_{X_2} &\lesssim& e^{a^\sharp t}\|(g,v)\|_\XX, 
\eean \Black
and the similar estimate for
$\|(S_{\BB  }\ast\AA 
S_{\BB })(t)(g,v)\|_\XX$. As a consequence, the announced estimate holds
for the family of operators $S_{\BB }\ast
\AA  S_{\BB }$.
\endproof

\subsection{The vanishing connectivity regime} 

When the network connectivity parameter vanishes, $\eps=0$, the linearized operator \Blue simplifies as \Black
\begin{equation}\label{LC1}
 {\varLambda_0}
    \begin{pmatrix}
    g \\ v
    \end{pmatrix}
    =
    \begin{pmatrix}
    -\partial_x g-a(x,0)g+\delta_{x=0}\OO _0[g,v]\\
    -\partial_y v  +\delta_{y=0}\OO _0[g,v]
    \end{pmatrix},
\end{equation}
where $\OO _0[g,v] = \NN_0[g] = \int_0^\infty a(x,0)g(x)\mathrm{d}x$. 
The {\Blue associated} semigroup is exponentially stable as shown in the following theorem.

\begin{theorem}\label{th:lc1}
There exist some constants $\alpha<0$ and $C>0$ such that
$\Sigma({\Blue\varLambda_0}) \cap \Delta_\alpha = \{ 0 \}$ and for any $(g_0, \, v_0) \in {\Blue\XX}$, $\langle g_0 \rangle = 0$, there holds 
\beqn\label{eq:lc}
\| S_{\Blue\varLambda_0}(t) (g_0, \,v_0) \|_{\Blue\XX} \leq C e^{\alpha t} \, \| (g_0, \, v_0) \|_{\Blue\XX}, \quad \forall \, t \ge 0.
\eeqn
\end{theorem}
\proof[Proof of Theorem~\ref{th:lc1}] 
Since $\Blue \varLambda^1_0=\Lambda_0$,  from Theorem \ref{lc1} we have already proved that  $g(t) := S_{\Blue\varLambda^1_0}(t)g_0$ satisfies 
$\|g(t)\|_{\Blue\XX_1}\leq C e^{{\Blue a} t}\,\|g_0\|_{\Blue\XX_1}$ for any
$t\geq0$ and any $a \in (a^*,0)$. We then  focus on
$\Blue\varLambda^2_0$. The Duhamel formula associated to the equation $\partial_t v=\Blue\varLambda^2_0 (g,v)$ writes
$$
v(t)=S_{\BB^2_0}(t) v_0+\int_0^t S_{\BB^2_0}(t-s)\AA^2_0\big(g(s), v(s)\big)\,\mathrm{d}s.
$$
Using the already known  estimate on $g(t)$, we deduce
\bean
\|S_{\Blue\varLambda^2_0}v_0(t)\|_{\Blue\XX_2} &=& \|v(t)\|_{\Blue\XX_2}\leq \|S_{\BB^2_0}(t)
v_0\|_{\Blue\XX_2}+\int_0^t \|S_{\BB^2_0}(t-s)\delta_0\NN_0[g(s)]\|_{\Blue\XX_2} \,
\mathrm{d}s\\
&\leq& e^{-\delta t} \|v_0\|_{\Blue\XX_2}+\int_0^t e^{-\delta(t-s)} C\,
e^{{\Blue a} s}\|g_0\|_{\Blue\XX_1}\, \mathrm{d}s\\
&\leq& C\, e^{\alpha t}\|(g_0,\,v_0)\|_{\Blue\XX}
\eean
for $\Blue\max\{a,-\delta\}<\alpha<0$, which yields our conclusion.
\endproof \Black

\subsection{Weak connectivity regime - exponential stability of the linearized equation}

In this part, we shall discuss the geometry structure of the
spectrum of the linearized time elapsed equation in weak
connectivity regime taking delay into account and using again a perturbation argument.

\begin{theorem}\label{th:CWR1}
There exists some constants $\eps_0 > 0$, $C \ge 1$ and $\alpha < 0$ such that for any $\eps\in[0,\eps_0]$
there holds $\Sigma({\color{red}\varLambda_\eps})\cap\Delta_\alpha=\{0\}$ and 
\beqn\label{Cwd}
\|S_{\color{red}\varLambda_\eps}(t)(g_0,v_0)\|_{\Blue\XX} \le C e^{\alpha t}\|(g_0,v_0)\|_{\Blue\XX},
\eeqn
for any $(g_0,v_0) \in {\Blue\XX}$ such that $\langle g_0 \rangle = 0$. 
\end{theorem}

\Blue We start presenting a technical result needed in the proof below. 
\Black
\begin{lemma}\label{lem:LCL1}
The operator $\Blue \varLambda_\eps$ is continuous with respect to
$\eps$, and more precisely
\begin{equation}\label{CoL1}
\|{\Blue\varLambda_\eps-\varLambda_0}\|_{\BBB(\Blue\XX)}\leq
O(\eps).
\end{equation}
\end{lemma}

\proof For all $(g,v)\in{\Blue\XX}$, we have
\begin{subequations}
 \begin{align}
 &{\Blue\varLambda_\eps}
 \begin{pmatrix}
 g \\ v
 \end{pmatrix}
 =
 \begin{pmatrix}
 -\partial_x g-a_\eps g-a'_\eps
 F_\eps\DD_\eps[v]+\delta_{x=0}\OO_\eps[g,v]\\
 -\partial_y v +\delta_{y=0}\OO_\eps[g,v]
 \end{pmatrix},\label{eq:Le1}\\
 &{\Blue\varLambda_0}
 \begin{pmatrix}
 g \\ v
 \end{pmatrix}
 =
 \begin{pmatrix}
 -\partial_x g-a(x,0) g +\delta_{x=0}\OO _0[g,v]\\
 -\partial_y v +\delta_{y=0}\OO _0[g,v]
 \end{pmatrix}. \label{eq:L01}
 \end{align}
\end{subequations}
From \eqref{eq:Le1}-\eqref{eq:L01}, we deduce
\begin{equation*}
({\Blue\varLambda_\eps-\varLambda_0})
\begin{pmatrix}
g \\ v
\end{pmatrix}
=
\begin{pmatrix}
(a(x,0)-a_\eps)g-a'_\eps
F_\eps\DD_\eps[v]+\delta_{x=0}(\OO_\eps[g,v]-\OO _0[g,v]) \\
\delta_{y=0}(\OO_\eps[g,v]-\OO _0[g,v])
\end{pmatrix}.
\end{equation*}
We then compute
\begin{eqnarray*}
\|({\Blue\varLambda_\eps-\varLambda_0})(g,v)\|_{\Blue\XX} &=&
\|(a(x,0)-a_\eps)g\|_{\Blue\XX_1}+\|a'_\eps
  F_\eps\DD_\eps[v]\|_{\Blue\XX_1}+2|\OO_\eps[g,v]-\OO _0[g,v]| \\
    &\leq& 3\|(a_\eps-a_0)g\|_{\Blue\XX_1}+2\|a'_\eps
    F_\eps\DD_\eps[v]\|_{\Blue\XX_2}\\
    &\leq& 3\eps\|a'\|_\infty\|g\|_{\Blue\XX_1}+2\eps a_1\|a'\|_\infty(1-\eps\|a'\|_\infty)\|F_\eps\|_{\Blue\XX_1} \|v\|_{\Blue\XX_2}\\
    &=& C\eps\|(g,v)\|_{\Blue\XX},
  \end{eqnarray*}
  which is nothing but \eqref{CoL1}.
\endproof
\proof[Proof of Theorem~\ref{th:CWR1}] With the help of Lemma~\ref{lem:LCL1}, we may proceed {\Blue similarly} as in the proof of Theorem~\ref{th:CWR} (see also again \cite{Tristani,MT*}) and we conclude that
$$
\Sigma({\Blue\varLambda_\eps})\cap\Delta_\alpha=\{\xi_\eps\},
$$
with $|\xi_\eps|\leq O(\eps)$ and $\xi_\eps$ is algebraically simple. We observe that 
$$
{\Blue\varLambda^*_\eps}
\begin{pmatrix}
\varphi \\ \psi
\end{pmatrix}
=
\begin{pmatrix}
\partial_x \varphi -a_\eps \varphi + a_\eps (\varphi(0) + \psi(0))
\\
\partial_y  \psi + \kappa_\eps b \, \psi(0) + \kappa_\eps b \, \varphi(0) - b \int a_\eps' F_\eps \, \varphi \, \mathrm{d}x 
\end{pmatrix},
$$
from which we deduce that ${\Blue\varLambda_\eps^*}(1,0) = 0$. Then  
$0\in\Sigma({\Blue\varLambda_\eps^*})$ and $\xi_\eps=0$. 
Moreover, the orthogonality condition  $\langle g_0 \rangle = \langle (g_0,v_0),(1,0) \rangle_{\Blue\XX,\XX'} =0$ implies that the exponential estimate \eqref{Cwd} holds.
%
\endproof

\subsection{Weak connectivity regime - nonlinear exponential stability}

We finally come back on the nonlinear problem and we present the proof of the second part of our main result for
the case with delay.

\proof[Proof of Theorem~\ref{th:MR} in case with delay] 
\
We write the system as
\bean
\partial_t f &=& - \partial_x f - a_\eps(\DD[u]) f + \delta_0 \PPP[f,\DD[u]]
 \\
 \partial_t u&=& - \partial_y u  + \delta_0 \PPP[f,\DD[u]], 
\eean
with 
$$
\PPP[f,m] = \int a(m) f, \quad \DD[u] = \int b u.
$$
We recall that the steady state $(F,U)$, $U := M {\bf 1}_{y \ge 0}$, satisfies
\bean
0 &=& - \partial_x F- a_\eps(M) F + \delta_0 M  
 \\
 0&=& - \partial_y U  + \delta_0 M, \quad M = \DD[U] = \PPP[F,\DD[U]].
\eean
We introduce the variation $g := f - F$ and $v = u - U$. The equation on $g$ is 
\bean
\partial_t g 
&=& - \partial_x g - a_\eps(\DD[u]) f + a_\eps(M) F + \delta_0 (\PPP[f,\DD[u]] - \PPP[F,\DD[U]])
\\ 
&=& - \partial_x g - a_\eps(M) f -a'_\eps F \DD[v] - Q[g,v] + \delta_0 \OO[g,v] + \delta_0 \QQ[g,v]
\\ 
&=& {\Blue\varLambda_\eps^1}(g,v) + \ZZ^1[g,v] , 
\eean
\Blue where
\bean
\QQ[g,v] &:=& \langle Q[g,v] \rangle,\\
\ZZ^1[g,v] &:=&  - Q[g,v] + \delta_0 \QQ[g,v],
\eean
with $Q[g,v]$ denoting that
\bean
Q[g,v]  &:=&  a_\eps(M) F - a_\eps(\DD[u]) f + a_\eps(M) f + a'_\eps F \DD[v]\\
 &=& \Phi(0) - \Phi(1) + \Phi'(0),
\eean
where
$$
\Phi({\Blue\theta})=a_\eps(\DD[{\Blue\theta}\,u+(1-{\Blue \theta})U])({\Blue\theta}\,f+(1-{\Blue\theta})F).
$$ 
\Black The equation on $v$ is 
\bean
\partial_t v 
&=& - \partial_y v  + \delta_0 (\PPP[f,\DD[u]] - \PPP[F,\DD[U]])\\ 
&=& - \partial_y v +  \delta_0 \OO[g,v] + \delta_0 \QQ[g,v]\\ 
&=& {\Blue\varLambda_\eps^2}(g,v) + \ZZ^2[g,v],
\eean
\Blue where
$$
\ZZ^2[g,v] := \delta_0 \QQ[g,v].
$$
\Black We then write the associated Duhamel formula
$$
(g(t),v(t)) = S_{\Blue\varLambda_\eps}(t) (g_0,v_0) + \int_0^t S_{\Blue\Lambda_\eps}(t-s) \ZZ[g(s),v(s)] \, \mathrm{d}s.
$$
{\Blue Refering to Step 2 in the proof of Theorem \ref{th:MR} gives a similar estimate $\|\ZZ[g,v] \|_{\Blue\XX} \le C\,\|(g,\,v)\|^2_{\Blue\XX}$, we then conclude the rest part of Theorem \ref{th:MR} in the case with delay}. 
\endproof

\bigskip 

\bibliographystyle{acm}

\begin{thebibliography}{}

\end{thebibliography}


\begin{thebibliography}{10}

\bibitem{MR0210343}
{\sc Coleman, B.~D., and Mizel, V.~J.}
\newblock Norms and semi-groups in the theory of fading memory.
\newblock {\em Arch. Rational Mech. Anal. 23\/} (1966), 87--123.

\bibitem{MR3311484}
{\sc De~Masi, A., Galves, A., L{\"o}cherbach, E., and Presutti, E.}
\newblock Hydrodynamic limit for interacting neurons.
\newblock {\em J. Stat. Phys. 158}, 4 (2015), 866--902.

\bibitem{MR1345150}
{\sc Diekmann, O., van Gils, S.~A., Verduyn~Lunel, S.~M., and Walther, H.-O.}
\newblock {\em Delay equations}, vol.~110 of {\em Applied Mathematical
  Sciences}.
\newblock Springer-Verlag, New York, 1995.
\newblock Functional, complex, and nonlinear analysis.

\bibitem{EngelNagel}
{\sc Engel, K.-J., and Nagel, R.}
\newblock {\em One-parameter semigroups for linear evolution equations},
  vol.~194 of {\em Graduate Texts in Mathematics}.
\newblock Springer-Verlag, New York, 2000.
\newblock With contributions by S. Brendle, M. Campiti, T. Hahn, G. Metafune,
  G. Nickel, D. Pallara, C. Perazzoli, A. Rhandi, S. Romanelli and R.
  Schnaubelt.

\bibitem{FL*}
{\sc Fournier, N., and L\"ocherbach, E.}
\newblock On a toy model of interacting neurons.
\newblock (2014) arXiv:1410.3263.

\bibitem{MR1554993}
{\sc Fredholm, I.}
\newblock Sur une classe d'\'equations fonctionnelles.
\newblock {\em Acta Math. 27}, 1 (1903), 365--390.

\bibitem{GKbook}
{\sc Gerstner, W., and Kistler, W.~M.}
\newblock {\em Spiking neuron models}.
\newblock Cambridge University Press, Cambridge, 2002.
\newblock Single neurons, populations, plasticity.

\bibitem{GMM*}
{\sc Gualdani, M.~P., Mischler, S., and Mouthot, C.}
\newblock Factorization of non-symmetric operators and exponential
  ${H}$-{T}heorem.
\newblock arXiv:1006.5523, 2010.

\bibitem{Kato}
{\sc Kato, T.}
\newblock {\em Perturbation theory for linear operators}.
\newblock Classics in Mathematics. Springer-Verlag, Berlin, 1995.
\newblock Reprint of the 1980 edition.

\bibitem{MS-erratum}
{\sc Mischler, S.}
\newblock Erratum: Spectral analysis of semigroups and growth-fragmentation
  equations.
\newblock submitted.

\bibitem{Mbook*}
{\sc Mischler, S.}
\newblock {\em Semigroups in Banach spaces - splitting approach for spectral
  analysis and asymptotics estimates}.
\newblock work in progress.

\bibitem{MMcmp}
{\sc Mischler, S., and Mouhot, C.}
\newblock Stability, convergence to self-similarity and elastic limit for the
  {B}oltzmann equation for inelastic hard spheres.
\newblock {\em Comm. Math. Phys. 288}, 2 (2009), 431--502.

\bibitem{MQT}
{\sc Mischler, S., Qui{\~n}inao, C., and Touboul, J.}
\newblock On a kinetic {F}itzhugh-{N}agumo model of neuronal network.
\newblock {\em Comm. Math. Phys. 342}, 3 (2016), 1001--1042.

\bibitem{MS}
{\sc Mischler, S., and Scher, J.}
\newblock Spectral analysis of semigroups and growth-fragmentation equations.
\newblock {\em Ann. Inst. H. Poincar\'e Anal. Non Lin\'eaire 33}, 3 (2016),
  849--898.

\bibitem{MT*}
{\sc Mischler, S., and Tristani, I.}
\newblock Uniform semigroup spectral analysis of the discrete, fractional \&
  classical fokker-planck equations.
\newblock (2015) hal-01177101.

\bibitem{PPS1}
{\sc Pakdaman, K., Perthame, B., and Salort, D.}
\newblock Dynamics of a structured neuron population.
\newblock {\em Nonlinearity 23}, 1 (2010), 55--75.

\bibitem{PPS2}
{\sc Pakdaman, K., Perthame, B., and Salort, D.}
\newblock Relaxation and self-sustained oscillations in the time elapsed neuron
  network model.
\newblock {\em SIAM J. Appl. Math. 73}, 3 (2013), 1260--1279.

\bibitem{PPS3}
{\sc Pakdaman, K., Perthame, B., and Salort, D.}
\newblock Adaptation and fatigue model for neuron networks and large time
  asymptotics in a nonlinear fragmentation equation.
\newblock {\em J. Math. Neurosci. 4\/} (2014), Art. 14, 26.

\bibitem{Pazy}
{\sc Pazy, A.}
\newblock {\em Semigroups of linear operators and applications to partial
  differential equations}, vol.~44 of {\em Applied Mathematical Sciences}.
\newblock Springer-Verlag, New York, 1983.

\bibitem{Q*}
{\sc Qui\~ninao, C.}
\newblock A microscopic spiking neuronal network for the age-structured model.
\newblock (2015) hal-01121061.

\bibitem{RibaricVidav}
{\sc Ribari{\v{c}}, M., and Vidav, I.}
\newblock Analytic properties of the inverse {$A(z)^{-1}$} of an analytic
  linear operator valued function {$A(z)$}.
\newblock {\em Arch. Rational Mech. Anal. 32\/} (1969), 298--310.

\bibitem{RT*}
{\sc Robert, P., and Touboul, J.~D.}
\newblock On the dynamics of random neuronal networks.
\newblock (2014) arXiv:1410.4072.

\bibitem{Tristani}
{\sc Tristani, I.}
\newblock Boltzmann equation for granular media with thermal force in a weakly
  inhomogeneous setting.
\newblock {\em J. Funct. Anal. 270}, 5 (2016), 1922--1970.

\bibitem{Voigt80}
{\sc Voigt, J.}
\newblock A perturbation theorem for the essential spectral radius of strongly
  continuous semigroups.
\newblock {\em Monatsh. Math. 90}, 2 (1980), 153--161.

\bibitem{MR0100765}
{\sc Volterra, V.}
\newblock {\em Theory of functionals and of integral and integro-differential
  equations}.
\newblock With a preface by G. C. Evans, a biography of Vito Volterra and a
  bibliography of his published works by E. Whittaker. Dover Publications,
  Inc., New York, 1959.

\bibitem{W*}
{\sc Weng, Q.}
\newblock General time elapsed neuron network model: well-posedness and strong
  connectivity regime.
\newblock (2015) hal-01243163.

\end{thebibliography}

\signsm  
\signqw

\end{document}